\documentclass[english,12pt,a4paper]{article}

\usepackage{xypic}
\usepackage{graphics}
\RequirePackage{amssymb} % for the square for the proof environment.
\RequirePackage{amsopn} % for DeclareMathOperator.
\RequirePackage{amsbsy} % for Boldsymbol
\RequirePackage{amsmath} % for multiple lines in equations
\RequirePackage{theorem} %
\RequirePackage{thc} %
\newcounter{mysubsection}[section]
\newcounter{mysubsubsection}[mysubsection]

\renewcommand{\thesection}%
                {\arabic{section}.}%
        \renewcommand{\thesubsection}%
                {\thesection\arabic{subsection}.}%
        \renewcommand{\thesubsubsection}%
                {\thesubsection\arabic{subsubsection}.}%
        \renewcommand{\themysubsection}%
                {(\thesection\arabic{mysubsection})}%
        \renewcommand{\themysubsubsection}%
                {(\thesection\arabic{mysubsection}.%
                \arabic{mysubsubsection})}%
\theoremstyle{change}%{changebreak}
\theoremheaderfont{\normalfont\bfseries}
\theorembodyfont{\upshape}

\newtheorem{remark}[mysubsection]{\remarkName}

\newtheorem{example}[mysubsection]{\exampleName}

\theorembodyfont{\slshape}
\newtheorem{theorem}[mysubsection]{\theoremName}
\newtheorem{proposition}[mysubsection]{\propositionName}
\newtheorem{lemma}[mysubsection]{\lemmaName}
\newtheorem{corollary}[mysubsection]{\corollaryName}
\def\qed{{\nobreak\hfil\penalty50% \unskip
  \hskip2em\hbox{}\nobreak\hfil$\square$%
  \parfillskip0pt\finalhyphendemerits=0\par}} % 3D0
\newenvironment{proof}%
  {\par\addvspace{\medskipamount}%
    \upshape%
    {\slshape\scshape%\bfseries%
    \proofName\hskip\labelsep}}%
  {\qed%
    \addvspace{\medskipamount}}%
\def\definitionName{Definition.}
\def\remarkName{Remark.}
\def\remarksName{Remarks.}
\def\noteName{Note.}
\def\notesName{Notes.}
\def\questionName{Question.}
\def\problemName{Problem.}
\def\exampleName{Example.}
\def\examplesName{Examples.}
\def\exerciseName{Exercise.}
\def\exercisesName{Exercises.}
\def\lpName{}
\def\hunchName{Hunch Hunch.}
\def\theoremName{Theorem.}
\def\propositionName{Proposition.}
\def\lemmaName{Lemma.}
\def\corollaryName{Corollary.}

\def\proofName{Proof.}

\DeclareMathOperator{\oF}       {\boldsymbol{\overline{F}}}%
\DeclareMathOperator{\f}        {\boldsymbol{f}}%
 \DeclareMathOperator{\Spec}   {Spec}
 \DeclareMathOperator{\Max}    {Max}
 \DeclareMathOperator{\Inv}{Inv}
 \DeclareMathOperator{\Na}{Na}
 \DeclareMathOperator{\Kr}{Kr}
 
 \DeclareMathOperator{\calB}   {\mathcal B}
 \DeclareMathOperator{\calM}   {\mathcal M}
 \DeclareMathOperator{\calW}   {\mathcal W}
 \DeclareMathOperator{\calT} {\mathcal T}
 \DeclareMathOperator{\calV} {\mathcal V}
\newcommand{\balf}
 {\renewcommand{\theenumi}{(\alph{enumi})}
 \renewcommand{\labelenumi}{\theenumi}
                      \begin{enumerate}}
\newcommand{\ealf}   {\end{enumerate}
                      \renewcommand{\theenumi}{\arabic{enumi}}
                      \renewcommand{\labelenumi}{\theenumi.}}
\newcommand{\bara}   {\renewcommand{\theenumi}{(\arabic{enumi})}
                      \renewcommand{\labelenumi}{\theenumi}
                      \begin{enumerate} }
\newcommand{\eara}   {\end{enumerate}
                      \renewcommand{\theenumi}{\arabic{enumi}}
                      \renewcommand{\labelenumi}{\theenumi.}}
\newcommand{\brom}   {\renewcommand{\theenumi}{(\roman{enumi})}
                      \renewcommand{\labelenumi}{\theenumi}
                      \begin{enumerate} }
\newcommand{\erom}   {\end{enumerate}
                      \renewcommand{\theenumi}{\arabic{enumi}}
                      \renewcommand{\labelenumi}{\theenumi.}}
\begin{document}

\title{Pr\"{u}fer $\star$--multiplication domains and \\ semistar operations}

\author{\textbf{M. Fontana}\thanks{Supported in part by a research grant MURST 2001/2002 (Cofin 2000-MM01192794).
 This project was partially developed during the visit of the first named author to
 Granada in the frame of the Socrates Program sponsored by EU.}\\
 Dipartimento di Matematica\\
 Universit{\`a} degli Studi Roma Tre\\
 Largo San Leonardo Murialdo, 1\\
 00146 Roma, Italy \\
 \textbf{P. Jara} and \textbf{E. Santos}\thanks{Partially supported by DGES BMF2001-2823 and FQM--266 (Junta de Andaluc{\'\i}a Research
 Group).}\\
 Department of Algebra\\
 University of Granada\\
 18071--Granada, Spain}

\date{}
\maketitle
\begin{abstract}
\noindent \footnotesize Starting from the notion of semistar operation, introduced in 1994 by
Okabe and Matsuda \cite{Okabe/Matsuda:1994}, which generalizes the classical concept of star
operation (cf.  Gilmer's book \cite{GILMER:1972}) and, hence, the related classical theory of
ideal systems based on the works by W. Krull, E. Noether, H. Pr\"{u}fer, P. Lorenzen and P. Jaffard
(cf.  Halter-Koch's book \cite{HK}), in this paper we outline a general approach to the theory of
Pr\"{u}fer $\star$-multiplication domains (or P$\star$MDs), where $\star$ is a semistar operation.
This approach leads to relax the classical restriction on the base domain, which is not
necessarily integrally closed in the semistar case, and to determine a semistar invariant
character for this important class of multiplicative domains (cf.  also J.M. Garc{\'\i}a, P. Jara and
E. Santos \cite{Garcia/Jara/Santos:1999}).  We give a characterization theorem of these domains
in terms of Kronecker function rings and Nagata rings associated naturally to the given semistar
operation, generalizing previous results by J. Arnold and J. Brewer \cite{Arnold/Brewer:1971} and
B.G. Kang \cite{Kang:1989}.  We prove a characterization of a P$\star$MD, when $\star$ is a
semistar operation, in terms of polynomials (by using the classical cha\-racte\-ri\-za\-tion of
Pr\"{u}fer domains, in terms of polynomials given by R. Gilmer and J. Hoffman
\cite{Gilmer/Hoffmann:1974}, as a model), extending a result proved in the star case by E.
Houston, S.J. Malik and J. Mott \cite{Houston/Malik/Mott:1984}.  We also deal with the
preservation of the P$\star$MD property by ``ascent'' and ``descent'' in case of field
extensions.  In this context, we generalize to the P$\star$MD case some classical results
concerning Pr\"{u}fer domains and P$v$MDs.  In particular, we reobtain as a particular case a result
due to H. Pr\"{u}fer \cite{Prufer} and W. Krull \cite{Krull:1936} (cf.  also F. Lucius
\cite{Lucius:1998} and F. Halter-Koch \cite{Koch:2000}). Finally, we develop several examples and
applications when $\star$ is a (semi)star given explicitly (e.g. we consider the case of the
``standard'' $v$--, $t$--, $b$--, $w$--operations or the case of semistar operations associated
to appropriate families of overrings).
\end{abstract}
\par
\textbf{Abbreviated title:} {Pr\"{u}fer $\star$--multiplication domains}
\section{Introduction}\label{intro}

The theory of ideal systems is based on the classical works by W. Krull, E. Noether, H. Pr\"{u}fer
and P. Lorenzen; a systematic treatment of this theory can be found in the volumes by P. Jaffard
\cite{Jaffard} and F. Halter-Koch \cite{HK}.  A different presentation, using the notion of star
operation, is given in 1972 by R. Gilmer \cite[Sections 32-34]{GILMER:1972} (cf. also for further
developments \cite {Hedstrom/Houston:1980}, \cite{Kang:1987}, \cite{Bachman:1988},
\cite{Anderson:1988}, \cite{Anderson/Anderson:1990}, \cite{OM1}, and \cite{Anderson/Cook:2000}).
In 1994 Okabe and Matsuda \cite{Okabe/Matsuda:1994} generalize the concept of star operation by
introducing the more ``flexible'' notion of semistar operation. After that paper new developments
of the multiplicative theory of ideals have been realized and successfully applied to analyze the
structure of different classes of integral domains (cf.  for instance \cite{Okabe/Matsuda:1997},
\cite{Matsuda/Sugatani:1995}, \cite{FH}, \cite{FL1}, \cite{FL2}, \cite{FL3}, \cite{DF:2001} and
\cite{Koch:JA}).

Semistar operations of a special type appear naturally in relation with the general constructions
of Kronecker function rings and Nagata function rings (in Section 1, we recall the definitions
and the principal properties of these objects).  More precisely, given a semistar operation
$\star$ on an integral domain $D$ with quotient field $K$, the Kronecker function ring $\Kr(D,
\star) \ (\subseteq K(X))\,$ [respectively, the Nagata function ring $\Na(D, \star) \ (\subseteq
K(X))\,$] induces naturally a ``distinguished'' semistar operation $\,\star_{a}\,$ [respectively,
$\,\tilde{\star}\,$] on $D$ such that $ \, F\Kr(D, \star) \cap K = F^{\star_{a}}\,$
[respectively, $\,\Na(D, \star) \cap K = F^{\tilde{\star}}\,$], for each finitely generated
fractional ideal $F$ of $D$.  These semistar operations were intensively studied in \cite{FL3},
where the authors examine also the interplay of $\Kr(D, \star)$ and $\star_{a}$ with
$\Na(D,\star)$ and $\tilde{\star}$ and show a ``parallel'' behaviour of these pairs of objects.

The equality of Nagata function ring with Kronecker function ring characterizes, in the classical
Noetherian case, the Dedekind domains.  It is natural, in the general context, to investigate on
the existence of ``semistar invariants'' for different classes of Pr\"{u}fer--like domains. A first
attempt in this direction is due to F. Halter-Koch \cite{Koch:2000}, who obtained a deep
axiomatic approach to the theory of Kronecker function rings, with applications to the
characterization of B\'ezout domains that are Kronecker function rings (cf.  also \cite{FL2}). On
the other hand, the study initiated in \cite{FL3} leads naturally  to the investigation of the
class of integral domains, having a semistar operation $\,\star\,$ such that the semistar
operation $\,\tilde{\star}\,,$ associated to the Nagata function ring, coincides with the
semistar operation $\,\star_{a}\,,$ associated to the Kronecker function ring.

One of the aims of this paper is to characterize a distinguished class of ``multiplication
domains'', called \it the Pr\"{u}fer semistar multiplication domains \rm or P$\star$MD, that
arises naturally in this context, having the property that
$\tilde{\star}=(\tilde{\star})_{a}=\star_a$ (Section 2).  This class contains as examples
Pr\"{u}fer domains, Krull domains and P$v$MD, but also integral domains, that are not integrally
closed, having although an appropriate overring which is Pr\"{u}fer star multiplicative domain
(cf. \cite{Houston/Malik/Mott:1984}, \cite{Kang:1989} and \cite{Garcia/Jara/Santos:1999}).  An
explicit example of a \it non \rm integrally closed Pr\"{u}fer semistar multiplication domain is
given in Example~\ref{ex:10} (recall that a Pr\"{u}fer star multiplication domain is always
integrally closed).

In Section 2 we show that, if $\,\star\,$ is semistar operation of finite type which is spectral
and e.a.b. on an integral domain $D$ (definitions are given in Section 1), then $D$ is a
P$\star$MD. Moreover we prove that $D$ is a P$\star$MD, for some semistar operation $\,\star\,$
on $D\,,$ if and only if $D$ is a P$\tilde{\star}$MD, where $\,\tilde{\star}\,$ is a semistar
operation of finite type which is spectral and e.a.b. This result extends one of the principal
results of \cite{Garcia/Jara/Santos:1999}, proved by using torsion theories.  After this
characterization, we apply our theory to some special types of semistar operations and we give
new characterizations of P$\star$MDs in the ``classical'' star setting.  In particular, we obtain
also that the P$w$MDs studied recently by W. Fanggui and R. L. McCasland
\cite{Fanggui/McCasland:1999} coincide with the P$v$MDs introduced by M. Griffin
\cite{Griffin:1967}.

In Section 3 we deal with the preservation of the
 P$\star$MD property by ``ascent'' and ``descent'', in case of algebraic field
 extensions.  We generalize to the P$\star$MD case some classical results
 concerning Pr\"{u}fer domains and P$v$MDs.  In particular, we reobtain
 the following generalization of a result due to H. Pr\"{u}fer  and W. Krull  (for the ``only if''
 case, cf. \cite[\S
 11]{Prufer} and \cite[Satz 9]{Krull:1936}) and to F. Lucius and F. Halter-Koch  (for
 the ``if'' case, cf.  \cite[Theorem 4.6 and Theorem
 4.4]{Lucius:1998} and \cite[Theorem 3.6]{Koch:2000}):
\newline
Let $K\subseteq{L}$ be an algebraic field extension, let
 $T$ be an integral domain with quotient field $L$, set $D:=T\cap{K}$.
 Assume that $D$ is integrally closed and that $T$ is the integral
 closure of $D$ in $L$.  Then $D$ is a P$v$MD if and only if $T$ is a
 P$v$MD.

\vskip0.1cm

We use as main reference Gilmer's book \cite{GILMER:1972} and any unexplained material is as in
\cite{GILMER:1972} and \cite{Kaplansky:1970}. Since many preliminary results on semistar
operations and applications, that we will need in this paper, are not easily available, because
the related work was presented or appeared in the Proceedings of recent Conferences (in
particular, \cite{FL1}, \cite{FL2} and \cite{FL3}), we will recall the principal definitions and
the statements of the main properties in Section 1.

\vskip0.1cm

The authors want to thanks Giampaolo Picozza for his very helpful comments and the
referee for his/her many valuable suggestions which have improved the previous version
of this paper.
\section{Background results}\label{sec:1}

Let $D$ be an integral domain with quotient field $K$. Let
$\boldsymbol{\overline{F}}(D)$ denote the set of all nonzero $D$-submodules of $K$ and let
$\boldsymbol{F}(D)$ be the set of all nonzero fractional ideals of $D$, i.e.  all $E \in
\boldsymbol{\overline{F}}(D) $ such that there exists a nonzero $d \in D$ with $dE \subseteq D$.
Let $\boldsymbol{f}(D)$ be the set of all nonzero finitely generated $D$-submodules of $K$.
Then, obviously $\, \boldsymbol{f}(D) \subseteq \boldsymbol{F}(D) \subseteq
\boldsymbol{\overline{F}}(D) \, .$

We recall that a mapping
\[
\star : \boldsymbol{\overline{F}}(D) \rightarrow \boldsymbol{\overline{F}}(D) \,, \;\;\;E \mapsto
E^{\star}  \] \noindent is called a \it semistar operation on $D$ \rm if, for $x \in K, x \not =
0$, and $E,F \in \boldsymbol{\overline{F}}(D) $, the following properties hold:

$(\star_1)$  $(xE)^{\star} = xE^{\star}\,; $

$(\star_2)$  $E \subseteq F \Rightarrow  E^{\star} \subseteq F^{\star}\,; $

$(\star_3)$  $E \subseteq E^{\star}$ and $E^{\star} = (E^{\star})^{\star}=:E^{\star \star} $

\noindent cf.  for instance \cite{OM1}, \cite{Okabe/Matsuda:1994}, \cite{Matsuda/Sugatani:1995},
\cite{Matsuda/Sato:1996}, \cite{FH} and \cite{FL1}.  In order to avoid trivial cases, we will
assume tacitly that the semistar operations are non trivial, i.e. if $D \not = K$ then $D^\star
\not = K$ (or, equivalently, the map $\star : \boldsymbol{\overline{F}}(D) \rightarrow
\boldsymbol{\overline{F}}(D)$ is not constant onto $K\,;\,$ cf.  \cite[Section 2]{FL3}).

A semistar operation $\star$ on $D$ is called an \it e.a.b. \rm (\it = endlich arithmetisch
brauchbar\rm) [respectively, \it a.b. \rm (\it = arithmetisch brauchbar\rm)] if, for each $E \in
\boldsymbol{f}(D)$ and for all $F,G \in \boldsymbol{f}(D)$ [respectively, $F,G \in
\boldsymbol{F}(D)]$: \[ (EF)^\star \subseteq (EG)^\star \Rightarrow F^\star \subseteq G^\star, \]
(cf.  for instance \cite[Definition 2.3 and Lemma 2.7]{FL1}).

If $\star_{1}$ and $\star_{2}$ are two semistar operation on $D$, we say that $\;
\star_{1} \leq \star_{2}\;$ if $E^{\star_{1}} \subseteq E^{\star_{2}}$, for each $E \in
\boldsymbol{\overline{F}}(D)\,;$ in fact, for semistar operations $\star_1$ and
$\star_2$, the following assertions are equivalent (i) $\star_1\leq\star_2$; (ii)
$(E^{\star_1})^{\star_2}=E^{\star_2}$ for each $E\in\oF(D)$ and (iii)
$(E^{\star_{2}})^{\star_{1}}=E^{\star_{2}}$ for each $E\in\oF(D)\,.$

Several new semistar operations can be derived
from a given semistar operation $\star$.  The essential details are given in the following
example.

\begin{example}\label{ex:1.1}
\rm  Let $D$ be an integral domain and let $\star$ be a semistar o\-pe\-ra\-tion on $D$.

\bf (a) \rm  If $\star$ is a semistar operation such that $D^\star = D\,$, then the map $\star :
\boldsymbol{F}(D) \rightarrow \boldsymbol{F}(D)\,$, \, $E \mapsto E^\star\,,$ is called a \it
star operation on \rm $D\,.$ Recall \cite[(32.1)]{GILMER:1972} that a star operation $\star$
verifies the properties $(\star_2)\,, (\star_3)\,,$ for all $E, F \in \boldsymbol{F}(D)\,;$
moreover, for each $x \in K\,, x \neq 0\,$ and for each $E \in \boldsymbol{F}(D)\,,$ a star
operation $\star$ verifies also:

\hspace*{10pt} $(\star\star_1)$ $(xD)^\star = xD\,, \; \; (xE)^{\star} = xE^{\star}\,.$

\noindent If $\star$ is a semistar operation on $D$ such that $D^\star = D\,,$ then we will write
often in the following of the paper that $\star$ is \it a (semi)star operation on $D\,,$ \rm for
emphasizing the fact that the semistar operation $\,\star\,$ is an extension to
$\boldsymbol{\overline{F}}(D)$ of \it a ``classical'' star operation \rm $\,\star\,,$ i.e. a map
$\, \star : \boldsymbol{F}(D)\rightarrow \boldsymbol{F}(D)\,,$ verifying the properties
$(\star\star_1)\,,(\star_2)\,$ and $\, (\star_3)\,$ \cite[Section 32]{GILMER:1972}.  Note that
not every semistar operation is an extension of a star operation \cite[Remark 1.5 (b)]{FH}.

\vskip.2cm
 \bf (b) \rm For each $E \in \boldsymbol{\overline{F}}(D)$, set
\[ E^{\star_f} := \cup \{F^{\star} \;|\;  \, F\subseteq E, \; F \in
\boldsymbol{f}(D) \}\,. \] \noindent Then $\star_f$ is also a semistar operation on $D$, called
\it the semistar operation of finite type associated to $\star\,$. \rm  Obviously, $F^{\star} =
F^{\star_{f}}$, for each   $F \in \boldsymbol{f}(D)\,.$ If $\star = \star_f$, then $\star$ is
called a \it  semistar operation of finite type \rm \cite[Example 2.5(4)]{FL1}. For instance, if
$\,v\,$ is \it the  $v$--(semi)star operation on $D$ \rm defined by $E^v := (E^{-1})^{-1},$ for
each $E \in \overline{\boldsymbol{F}}(D)\,,$ with $E^{-1} := (D :_{\mbox{\tiny \it K}} E) := \{ z
\in K \; | \;\; zE \subseteq D \}\,)\,$ \cite[Example 1.3 (c) and Proposition 1.6 (5)]{FH},\,
then the semistar operation of finite type $\,v_{f}\,$ associated to $\,v\,$ is called \it the
$t$--(semi)star operation on $D$\, \rm (in this case $D^v =D^t = D$).
\newline
Note that, in general,  $\star_f \leq \star\,,$ i.e. $\, E^{\star_f} \subseteq E^{\star}\, $ for
each $E \in \boldsymbol{\overline{F}}(D)$.  Thus, in particular, if $E = E^{\star}$, then $E =
E^{\star_f}$.  Note also that $\star_f = (\star_f)_{f}$.
\newline
We say that \it two semistar operation on \rm $D$, $\star_{1}$ and $\star_{2}$, are \it
equivalent \rm if $(\star_{1})_{f} = (\star_{2})_{f}\,.$

\vskip.2cm
 \bf (c) \rm Next example of a semistar operation is connected with the constructions
already in \cite{Zafrullah:1985}, \cite{Anderson:1988} and \cite{Anderson/Anderson:1990} and with
a weak version of integrality.  The essential techniques and motivations for considering this
weak version of integrality, using ideal systems, can be found in Jaffard's book \cite{Jaffard}.
More recently, starting from an idea in \cite{Anderson/Houston/Zafrullah:1991}, where the authors
introduced a weak version of integrality (called semi-integrality and associated to the
$v$--operation), a weak general version of integrality, depending on a star operation, was
introduced and studied in \cite{OM1}, \cite{Koch:1997}, \cite{Koch:2000} and
\cite{Fanggui:2001}.  The natural extension of this notion to the case of semistar operations was
considered in \cite{Okabe/Matsuda:1994}, \cite{FL1} and \cite{FL2}.
\newline
We start by defining a new operation on $D$, denoted by $[\star]$, called \it the semistar
integral closure of $\star$,\rm $\,$ by setting:
$$
 F^{[\star]} := \cup\{((H^\star:H)F)^{\star_f} \,\;|\;\, H \in \boldsymbol{f}(D) \} \,,\;\;
 \mbox{for each } \, %^\star$
 F \in \boldsymbol{f}(D) \,,
$$
and
$$
 E^{[\star]} := \cup \{F^{[\star]} \,\;|\;\, F \in
 \boldsymbol{f}(D),\;
  F \subseteq E\}\,,\; \; \mbox{for each } \, E \in
 \boldsymbol{\overline{F}}(D) \,.
$$
It is not difficult to see that the operation $[\star]$
defined in this manner is a semistar operation of finite type on $D$, that $\, \star_{f} \leq
[\star]\,$, hence $\, D^{\star} \subseteq D^{[\star]}\,$, and that $D^{[\star]}$ is integrally
closed \cite[Definition 4.2, Proposition 4.3 and Proposition 4.5 (3)]{FL1}.  Therefore, it is
obvious that if $\, D^{\star} = D^{[\star]}\, $ then $D^{\star}$ is integrally closed.  The
converse is false, even when $\, \star\,$ is a (semi)star operation on $D$.

\bf (c.1) \rm There exists \it an integral domain $D$ with a semistar operation $\star$ such that
$D^{\star}$ is integrally closed and $D^{\star} \subsetneq D^{[\star]}$.  \rm
\newline
Let $V$ be a valuation domain of the form $\,K+M\,$, where $K$ is a field and $M$ is the maximal
ideal of $V$.  Let $k$ be a proper subfield of $K$ and assume that $k$ is algebraically closed in
$K$.  Set $D:= k+M \subsetneq V$ and consider the (semi)star operation $\, \star := v \, $ on
$D$.  Then, clearly, $D$ is integrally closed and $D^\star (= D^{v}) = D$.  On the other hand,
let $\, z \in K \setminus k\,$ and let $W:=k+zk$ then $W$ is a $k$--submodule of $K$, which
obviously is not a fractional ideal of $k$. Then $\,H:=W+M\,$ is a finitely generated fractional
ideal of $D$ and $H^{v }= V$ by \cite[Theorem 4.3 and its proof]{Bastida/Gilmer:1974}. Therefore
$(H^{v}:H^{v }) =V$, and so $ V \subseteq D^{[v]}$ (in fact, $ V = D^{[v]}$ by \cite[Proposition
8 (ii)]{Anderson/Houston/Zafrullah:1991}).
\newline
A simple case for having that $D^{\star}$ is integrally closed if and only if $\, D^{\star} =
D^{[\star]}\,$ is when $\,\star\,$ is a semistar operation of finite type on $D$ which is stable
with respect to finite intersections (i.e. $\, (E \cap F)^{\star} = E^{\star} \cap F^{\star}\,,$
for all $E, F \in \boldsymbol{\overline{F}}(D)\,)$.

\bf (c.2) \it  Let $\,\star\,$ be a semistar operation of an integral domain $D$.
Assume that $\star_{f}$ is stable, then $D^{[\star]}=(D')^{\star_f}$, where $D'$ is the
integral closure of $D$.\rm
\newline
Indeed,
 $D^{[\star]}
 =\cup\{(H^{\star_f}:H)\mid\;H\in\f(D)\}
 =\cup\{(H:H)^{\star_f}\mid\;H\in\f(D)\}
 \subseteq(D')^{\star_f}
 \subseteq(D^{[\star]})^{[\star]}=D^{[\star]}$.
\newline
In particular, if $D^\star$ is integrally closed, then
$D\subseteq{D'}\subseteq{D^\star}$  implies $D^{[\star]}=(D')^{\star_f}=D^\star$.

\vskip.2cm
 \bf (d) \rm The essential constructions related to the following example of
semistar operation  are due to P. Lorenzen \cite{Lorenzen} and P. Jaffard \cite{Jaffard} (cf.
also F. Halter-Koch \cite{HK}).
\newline
Given an arbitrary semistar operation $\star$ on an
 integral domain $D$, it is possible to associate to $\star$, an e.a.b.
 semistar operation of finite type $\, \star_a \, $ on $D$, called \it
 the e.a.b. semistar
operation associated to $\star$, \rm defined as follows:
$$
 F^{\star_a} :=
 \cup\{((FH)^\star:H)\,\; |\; \, H \in \boldsymbol{f}(D) \}, \; \; \textrm {for each } \, %^\star
 F \in \boldsymbol{f}(D) \, ,
$$
and $$ E^{\star_a} := \cup\{F^{\star_a} \,\; |\; \, F \subseteq
E\,,\; F \in \boldsymbol{f}(D) \},\; \; \textrm {for each } \,  E \in
\boldsymbol{\overline{F}}(D), $$ \cite [Definition 4.4]{FL1}.  Note that $\, [\star] \leq
{\star_a}\,$, that $D^{[\star]} = D^{\star_a}$ and if $\star$ is an e.a.b. semistar operation of
finite type then $\,\star = \star_{a}\,$ \cite [Proposition 4.5]{FL1}.

\vskip.2cm
\bf (e) \rm   Let $D$ be an integral domain and $T$ an overring of $D$.  Let
$\star$ be a semistar operation on $D$ and define
  $\dot{\star}^{\mbox{\tiny \it \tiny T}}: \overline{\boldsymbol{F}}(T)
  \rightarrow \overline{\boldsymbol{F}}(T)$ by setting:
$$ E^{\dot{\star}^{\mbox{\tiny \it \tiny T}}} := E^\star\,, \; \mbox{ for each } \;  E \in
\overline{\boldsymbol{F}}(T) (\subseteq \overline{\boldsymbol{F}}(D))\,. $$ Then, we know
\cite[Proposition 2.8]{FL1}:

\bf (e.1) \it The operation $\dot{\star}^{\mbox{\tiny \it T}}$ is a semistar operation on $T$
and, if $\star$ is of finite type on $D$, then ${\dot{\star}}^{\mbox{\tiny \it T}}$ is also of
finite type on $T\,.$

\bf (e.2) \it When $T = D^\star$, then ${\dot{\star}}^{{\mbox{\tiny \it D}}^{\star}}$, restricted
to $\boldsymbol{F}(D^\star)\,,$ defines a star operation on $D^\star\,.$

\bf (e.3) \it If $\star$ is e.a.b., then ${\dot{\star}}^{{\mbox{\tiny \it \tiny D}}^{\star}}$ is
also e.a.b.

\rm Conversely, let $\star$ be a semistar operation on an overring $T$ of $D$ and
define \ \d{$\star$}$_{\mbox{\tiny \it \tiny D}}: \overline{\boldsymbol{F}}(D) \rightarrow
\overline{\boldsymbol{F}}(D)$ by setting: $$ E^{{\mbox{\d{$\star$}}}_{\mbox{\tiny \it \tiny D}}}
  := (ET)^\star\,, \; \mbox{ for each } \; E \in
\overline{\boldsymbol{F}}(D) \,. $$ Then, we know \cite[Proposition 2.9, Corollary 2.10]{FL1}:

\bf (e.4) \it The operation \d{$\star$}$_{\mbox{\tiny \it \tiny D}}$ is a semistar operation on
$D\,.$

\bf (e.5) \it If $\star$ is e.a.b., then \d{$\star$}$_{\mbox{\tiny \it \tiny D}}$ is also e.a.b.

\bf (e.6) \it If we denote simply by $\boldsymbol{\ast}$ the semistar operation
\d{$\boldsymbol{\star}$}$_{{\mbox{\tiny \it \tiny D}}}\,$, then the semistar operations \,
$\dot{\boldsymbol{\ast}}^{\mbox{\tiny \it T}}\,$ and $\, \star \, $ (both defined on $T\,$)
coincide.

\noindent \rm Note that the module systems approach, developed by Halter-Koch in
\cite{Koch:JA}, gives a natural and general setting for (re)considering semistar
operations and, in particular, the semistar operations $\dot{\star}^{\mbox{\tiny \it
T}}$ and \d{$\star$}$_{\mbox{\tiny \it \tiny D}}$ .

\vskip.2cm
\bf (f) \rm Let $\Delta $ be a nonempty set of prime ideals of an integral domain $D$.  For each
$D$-submodule  $E$ of $K\,,$ set: $$E^{\star_\Delta}:=\cap \{ ED_P\; | \;\, P\in \Delta \}\,.$$
The mapping $E\mapsto E^{\star_\Delta}$, for each $E\in \boldsymbol{\overline{F}}(D)$, defines a
semistar operation on $D\,,$ moreover \cite [Lemma 4.1]{FH}:

\bf (f.1) \it For each $E\in \boldsymbol{\overline{F}}(D)$ and for each $P\in \Delta\,$,
$ED_P=E^{\star_\Delta}D_P$.

\bf (f.2) \it The semistar operation $\star_\Delta$ is \it stable (with respect to the finite
intersections), \it i.e. for all $E, F \in \boldsymbol{\overline{F}}(D)\,$ we have $\, (E \cap
F)^{\star_{\Delta}} = E^{\star_{\Delta}} \cap F^{\star_{\Delta}}\,.  $

\bf (f.3) \it For each $P\in \Delta$, $P^{\star_\Delta}\cap D = P$.

\bf (f.4) \it For each nonzero integral ideal $I$ of $D$ such that $I^{\star_\Delta}\cap D\neq
D$, there exists a prime ideal $P\in \Delta$ such that $I\subseteq P$.

\vskip.2cm
A \it semistar operation $\star$ \rm is called \it spectral, \rm if
there exists a nonempty set $\Delta $ of Spec($D$) such that $\star = \star_{\Delta}\,;$  in this
case we say that $\star$ is \it the spectral semistar operation associated with \rm $\Delta \,.$
We say that $\star$ is a \it quasi--spectral semistar operation \rm (or that \it $\star$ possesses
enough primes\ \rm) if, for each nonzero integral ideal $I$ of $D$ such that
$I^{\star_\Delta}\cap D\neq D$, there exists a prime ideal $P$ of $D$ such that $I\subseteq P$
and $P^{\star}\cap D = P\,.$ From (f.3) and (f.4), we deduce that each spectral semistar
operation is quasi--spectral.

A subset $\Delta $ of Spec($D$) is called \it stable for generizations \rm if $ \
Q \in \Spec(D)\ ,$ \; $ P \in \Delta \, $ and $\, Q \subseteq P\,, $ then $\, Q \in \Delta\,.$
Set $ \Delta^{\downarrow}:= \{ Q \in \Spec(D) \; | \;\; Q \subseteq P\, \mbox{ for some } P \in
\Delta \}\,$ and let $\Lambda \subseteq \Spec(D),$ it is easy to see that:

\bf (f.5) \it If $\, \Delta \subseteq \Lambda \subseteq \Delta^{\downarrow}\,,$ \, then
\;$\star_{\Delta} = \star_{\Lambda} = \star_{\Delta^{\downarrow}}\,.$

\vskip.2cm
\bf (g) \rm Example (f) can be generalized as follows. Let $ \calT := \{ T_\alpha
\; | \;\; \alpha \in A \}$ be a nonempty family of overrings of $D\,$ and define $\ \star_{\calT}:
\overline{\boldsymbol{F}}(D) \rightarrow \overline{\boldsymbol{F}}(D)$ by setting: $$
E^{\star_{\calT}} := \cap \{ ET_\alpha \; | \;\;  \alpha \in A \} \,,
  \; \mbox{ for each } \;  E \in  \overline{\boldsymbol{F}}(D) \,.
$$ Then we know that \cite[ Lemma 2.4 (3), Example 2.5 (6), Corollary 3.8]{FL1}:

\bf (g.1) \it The operation $ \star_{ {\calT}} $ is a semistar operation on $D$. Moreover, if $
{\calT}= \{ D_P \; | \;\; P \in \Delta \},$ then $ \star_{\calT} = \star_{\Delta}\,.$

\bf (g.2) \it  $ E^{\star_{\calT} }T_\alpha = E T_\alpha \,, \;$ for each $E \in
\overline{\boldsymbol{F}}(D)\,$ and for each $\alpha \in A\,.$

\bf (g.3) \it  If $ {\calT} = {\calW}$ is a family of valuation overrings of $D$, then
$\star_{\calW}$ is an a.b.  semistar operation on $D$.  If  ${\calW}$ is the family of \it all
\rm the valuation overrings of $D$, then $\star_{\calW}$ is called \it the $b$--semistar
operation on $D\,;\;$\rm  moreover, if $D$ is integrally closed, then $D^b = D\,$ \cite[Theorem
19.8]{GILMER:1972}, and thus the operation $b\,,$ restricted to $ \boldsymbol{F}(D)\,,$ defines a
star operation on $D\,,$ called \it the $b$-star operation \rm \cite[p. 398]{GILMER:1972}.  \rm

\vskip 0.2cm
Let $\star$ be a semistar operation of an integral domain $D$ and
assume that the set: $$
  \Pi^\star := \{ P \in \mbox{Spec}(D) \; | \;\; P\neq 0 \mbox{ and }
  P^{\star} \cap D \neq D \}
  $$
is nonempty, then the spectral semistar operation of $D$ defined by $\, \star_{sp}:=
\star_{\Pi^{\star}}\, $ is called \it the spectral semistar operation associated to $\star\,.$
\rm \, Note that if $\star$ is quasi--spectral, then $ \Pi^{\star}$ is nonempty and $ \star_{sp}
\leq \star \, $ \cite[Proposition 4.8 and Remark 4.9]{FH}.  It is easy to see that $\, \star \,$
is spectral if and only if $\, \star = \star_{sp}\,.$
\end{example}
\medskip

Let $I\subseteq D$ be a nonzero ideal of $D$.  We say that $I$ is a \it
quasi--$\star$--ideal of $D$\,  \rm if $I^\star \cap D = I\,.\;$
Note that, for each nonzero integral ideal $I$ of $D\,,$ the
  ideal $J:= I^\star \cap D$ is a quasi--$\star$--ideal of $D\,$ and $\,I
  \subseteq J \,.$
 Note also that the quasi--$\star$--ideals form a weak ideal system on
  $D$, in the sense of \cite{HK}: this alternative approach can be applied
  for recovering some of the results mentioned next.

  \noindent A \it quasi--$\star$--prime \ \rm [respectively, a \it
  quasi--$\star$--maximal\ \rm] is a quasi--$\star$--ideal which is also a
  prime ideal [respectively, quasi--$\star$--ideal which is a maximal
  element in the set of all proper quasi--$\star$--ideals of $D\,$].  \,
  It is not difficult to see that,

  %LEMMA 1.2
  \begin{lemma} \label{lm:1.2} {\rm \cite [Lemma 2.4]{FL3}}
When $\star =\star_f$, then:

{\bf (a)} each proper quasi--$\star$--ideal is contained in a quasi--$\star$--maximal;

{\bf (b)} each quasi--$\star$--maximal is a quasi--$\star$--prime;

{\bf (c)} the (nonempty) set $\calM(\star)$ of all quasi--$\star$--maximals coincide with the set:
   $$\Max\{P\in \Spec(D) \; |
   \;\; 0\neq{P}\mbox{ \it and } P^{\star} \cap{D} \not = D\} =
   \Max(\Pi^{\star})\,.$$ \hfill $\square$
\end{lemma}

  %REMARK 1.3
     \begin{remark} \label{rk:1.3} \rm Note that, if $\star$ is a
semistar
     operation of finite type, then $\star$ is quasi--spectral
     (Lemma~\ref{lm:1.2} ((a) and (b))).  Moreover, by Lemma~\ref{lm:1.2} (c) and
Example~\ref{ex:1.1} (f.5),
\[
 (\star_{f})_{sp } = \star_{\calM(\star_{f})}\,.
\]
We will simply denote by $\, \tilde{\star} \,$ the spectral semistar  operation
$\,(\star_{f})_{sp}\,,$ (cf.  also \cite[Proposition 3.6 (b) and Proposition 4.23 (1)]{FH}). From
the previous considerations it follows that $\, \tilde{\star} \leq \star_{f} \,$ and that
$\,\tilde{\star}\,$ is a spectral semistar operation of finite type (cf. also
\cite[Proposition~3.2 and Theorem 4.12 (2)]{FH}).
\newline
When $\,\star\,$ is the (semi)star
     $v$--operation, the (semi)star  operation $\,\tilde{v}\,$ coincides with
     \it the (semi)star operation $\,w\,$ \rm defined as follows: $$E^{w} :=
     \cup \{(E:H)\; | \;\,H \in \boldsymbol{f}(D) \mbox{ and } H^v = D \}\,,
     \;\; \mbox{ for each } E \in \boldsymbol{\overline{F}}(D)\,.$$
     This (semi)star operation was firstly considered by J. Hedstrom and E.
     Houston in 1980 \cite[Section 3]{Hedstrom/Houston:1980} under the name
     of F$_{\infty}$--operation, starting from the F--operation introduced by
     H. Adams \cite{Adams}.  Later, from 1997, this operation was
     intensively studied by W. Fanggui and R. McCasland (cf.
     \cite{Fanggui/McCasland:1997}, \cite{Fanggui/McCasland:1999} and
     \cite{Fanggui:2001}) under the name of $w$--operation.  Note also that
     the notion of $w$--ideal coincides with the notion of semi-divisorial
     ideal considered by S. Glaz and W. Vasconcelos in 1977
     \cite{Glaz/Vasconcelos:1977}.  Finally, in 2000, for each (semi)star
     operation $\,\star\,$, D.D. Anderson and S.J. Cook
     \cite{Anderson/Cook:2000} considered the $\,\star_{w}$--operation which
     can be defined as follows: $$E^{\star_{w}} := \cup \{(E:H)\; | \;\,H \in
     \boldsymbol{f}(D) \mbox{ and } H^\star = D \}\,,
     \;\; \mbox{ for each } E \in \boldsymbol{\overline{F}}(D)\,.$$
     From their theory it follows that $\,\star_{w} = \tilde{\star}$\,
     \cite[Corollary 2.10]{Anderson/Cook:2000}.  A deep link between the
     semistar operations of type $\,\tilde{\star}\,$ and the localizing systems
     of ideals was established in \cite{FH}.

    \end{remark}

\vskip 0.3cm
Let $R$ be a ring and $X$ an indeterminate over $R\,,$ for each $f
\in R[X]\,,$ we denote by $\boldsymbol{c}(f)$ the \it content of \rm $f$, i.e. the ideal of $R$
generated by the coefficients of the polynomial $f$.  The following ring, subring of the total
ring of rational functions: $$R(X) := \left\{ \frac{f}{g} \; | \;\; f, g \in R[X] \, \mbox{ and }
\, \boldsymbol{c}(g) = R \right\}$$
    is called the \it Nagata ring of \rm $R\,$
     \cite[Proposition 33.1]{GILMER:1972}.

     %LEMMA 1.4

    \begin{lemma} \label{lm:1.3} \rm \cite[Proposition 3.1]{FL3} \it
Let $\star$ be a semistar operation of an integral domain $D$ and set: $$
    N(\star) := N_{D}(\star) := \{ h \in D[X] \; | \;\; h \not = 0
\mbox{
    and } \,
    \boldsymbol{c}(h)^\star = D^\star \}\,.$$

\bf (a) \it $N(\star) = D[X] \setminus \left(\cup\{Q[X]
    \; | \;\; Q \in \calM(\star_{f}) \}\right)$ is a saturated
multiplicatively
    closed subset of $D[X]$ and, obviously, $N(\star) = N(\star_{f})$.

\bf (b) \it $\Max(D[X]_{N(\star)}) = \{Q[X]_{N(\star)} \; | \;\; Q \in \calM(\star_{f}) \}$.

\bf (c) \it $D[X]_{N(\star)} = \cap\{ D[X]_{Q[X]} \; |
    \;\; Q \in \calM(\star_{f}) \} = \cap\{ D_{Q}(X) \; | \;\; Q \in
    \calM(\star_{f}) \}$.

\bf (d) \it $\calM(\star_{f})$ coincides with the canonical image into \rm Spec($D$) \it of the
set of the maximal ideals of $D[X]_{N(\star)}\,,$ i.e. $\calM(\star_{f}) = \{ M \cap D \; | \; \;
M \in \Max({D[X]}_{N(\star)})\}\,).$\hfill $\square$
    \end{lemma}

We set:
    $$ \mbox{Na}(D, \star) := D[X]_{N_{D}(\star)}$$
    and we call this integral domain \it  the Nagata ring of $D$ with
respect to the semistar operation $\,\star$ \rm.  Obviously, $ \mbox{Na}(D,\star) = \mbox{Na}(D,
\star_{f})$ and if $\,\star = d\,,$ where $d$ is \it the identical (semi)star operation of \rm
$D$ (i.e. $E^d := E\,,$ \, for each $E \in \boldsymbol{\overline{F}}(D)),$ \ then $ \mbox{Na}(D,
d) = D(X)\,.$

%LEMMA 1.5
\begin{lemma} \label{lm:1.4} \rm \cite[Corollary 2.11,
Proposition 3.4, Corollary 3.6, Theorem 3.9]{FL3} \it Let $\,\star\,$ be a given semistar
operation of an integral domain $D$ and let $\,\tilde{\star} := \star_{\calM(\star_{f})} =
(\star_{f})_{sp}\,$ be the spectral semistar operation of finite type canonically associated to
$\star$ (cf. Remark~\ref{rk:1.3})).  Denote simply by $\dot{\tilde{\star}}$ the following
(semi)star operation on $D^{\tilde{\star}}$ (Example~\ref{ex:1.1} (e)): $$
{\dot{\tilde{\star}}}^{\mbox{{\tiny D}}^{\tilde{\star}}}:
\overline{\boldsymbol{F}}(D^{\tilde{\star}}) \rightarrow
\overline{\boldsymbol{F}}(D^{\tilde{\star}}) \,, \;\; E \mapsto E^{\tilde{\star}}\,. $$
 Then, for each $E \in \boldsymbol{\overline{F}}(D)$,

\bf (a) \it $E^{\star_{f}} = \cap \{ E^{\star_{f}}D_{Q}\; | \;\; Q \in \calM(\star_{f}) \}\,;$

\bf (b) \it  $ E^{\tilde{\star}} = \cap \{ ED_{Q}\; | \;\; Q \in \calM(\star_{f}) \}\,;$

\bf (c) \it $E\mbox{\rm Na}(D, \star) = \cap \{ ED_{Q}(X)\; | \;\; Q \in \calM(\star_{f}) \}\,;$

\bf (d) \it $E\mbox{\rm Na}(D, \star) \cap K = \cap \{ ED_{Q}\; | \;\; Q \in \calM(\star_{f})
\}\,;$

\bf (e) \it $E^{\tilde{\star}} = E\mbox{\rm Na}(D, \star) \cap K\,.$

\bf (f) \it For each $\ Q \in \calM(\star_{f})\,,$ set $\, Q^\diamond := QD_{Q}(X) \cap \Na(D,
\star)\,,$ then $\, Q^\diamond = Q[X]_{N_{D}(\star)}\in \Max(\Na(D, \star))$\, and $\,
\Na(D,\star)_{Q^\diamond} = D_{Q}(X)\,.$

\bf (g) \it $\calM(\star_{f}) = \calM(\tilde{\star}).$

\bf (h) \it $\calM(\dot{\tilde{\star}}) = \{\widetilde{Q} := QD_{Q}\cap D^{\tilde{\star}} \; |
\;\; Q \in \calM(\star_{f})\}\,$ and $\, D_{\widetilde{Q}}^{\tilde{\star}} = D_{Q}\,,$ for each
$\ Q \in \calM(\star_{f})\,.$

\bf (i) \it $\Na(D, \star) = \Na(D, {\tilde{\star}}) = \Na(D^{\tilde{\star}},\dot{\tilde{\star}})
\supseteq D^{\tilde{\star}}(X)\,.$ \hfill $\square$
\end{lemma}

\smallskip

We recall now a notion of invertibility that generalizes the classical concepts of invertibility,
$v$--invertibility and $t$--invertibility (cf.  for instance \cite{Anderson/Zafrullah:1997} and
\cite[Section 2]{Anderson/Cook:2000}).  Let $\star$ be a semistar operation on an integral domain
$D$.  Let $I \in \boldsymbol{F}(D)$, we say $I$ is $\star$--{\em invertible} if $(II^{-1})^\star
= D^\star$. Note that, if $I\in \boldsymbol{f}(D)$, then $ I\mbox{ is }\star_f\mbox{--invertible
}$ if and only if there exists $ J \in \boldsymbol{f}(D) \mbox{ such that }(IJ)^\star= D^\star$
and $J \subseteq I^{-1}\,,$ \, \cite{DF:2001}.  The following lemma ge\-ne\-ra\-lizes a result
proved by B.G. Kang \cite[Theorem 2.12]{Kang:1989} (cf.  also
 D.D. Anderson \cite[Theorem 2]{Anderson:1976}).

%LEMMA 1.6
\begin{lemma}\label{le:1} \rm \cite[Theorem 2.5]{DF:2001} \it Let $\star$ be a
semistar operation on an integral domain $D$.  Assume that $\star = \star_f$.  Let $I \in
\boldsymbol{f}(D) $, then the following are equivalent:
\brom \rm
\item \it $I$ is $\star$--invertible;

\rm \item \it $ID_Q\in\Inv(D_Q)$, for each $Q\in\calM(\star)$;

\rm \item \it $I\Na(D,\star)\in\Inv(\Na(D,\star))$.
\erom
\vskip -0.5cm \hfill $\square$
\end{lemma}

Let $\star$ be a semistar operation on an integral domain $D$.  We say that $D$ is a P$\star$MD
({\em Pr\"{u}fer $\star$--multiplication domain}), if each $I\in\boldsymbol{f}(D)$ is
$\star_f$--invertible.
\newline
It is obvious that if $\,\star_{1} \leq \star_{2}\,$ are two semistar operations on
an integral domain $D$ and if $D$ is a P$\star_{1}$MD, then $D$ is also a P$\star_{2}$MD.
Moreover, if $\star_{1}$ is equivalent to $\star_{2}\,,$ then $D$ is a P$\star_{1}$MD if and only
if $D$ is also a P$\star_{2}$MD. In particular, the notions of P$\star$MD and P$\star_{f}$MD
coincide.
\newline
Note that if $\star$ is a semistar operation on $D$ such that $D^\star = D$ (i.e. if $\star$
restricted to $\boldsymbol{F}(D)$ defines a star operation on $D\,;$ cf. Example~\ref{ex:1.1}
(a)), then $\star \leq v\,$ (where $v$ is the $v$--(semi)star operation, Example~\ref{ex:1.1}
(b)) \cite[Theorem 34.1 (4)]{GILMER:1972}. In particular, if $D^\star = D\,,$ then $\star_{f}\leq
t\,$ (where $t$ is the (semi)star operation of finite type associated to $v$); moreover, in the
present situation, if $D$ is a P$\star$MD, then $D$ is also a P$v$MD. In the semistar case a
P$\star$MD is not necessary a P$v$MD (see Example~\ref{ex:10} below).
\newline
Recall that if $d$ is the identical (semi)star operation on $D$, then obviously $d \leq \star$,
for each semistar operation $\star$ on $D\,.$ \, Moreover, the notion of P$d$MD coincide with the
notion of a Pr\"{u}fer domain \cite[Theorem 21.1]{GILMER:1972}.  Therefore, a Pr\"{u}fer domain is a
P$\star$MD, for each semistar operation $\star\,.$

%LEMMA 1.7
  \begin{lemma}\label{le:K} \rm (\cite[Theorem 3.11 (2), Theorem 5.1,
  Corollary 5.2, Corollary 5.3]{FL1}) \it Let $\star$ be any semistar operation
  defined on an integral domain $D$ with quotient field $K$ and let
  $\star_{a}$ be the e.a.b. semistar operation associated to $\star\,$
  (Example~\ref{ex:1.1} (d)).  Consider the e.a.b. (semi)star operation
  $\dot{\star}_{a} :={{\dot{\star}}_{a}}^{{\mbox{\tiny \it \tiny
  D}}^{\star_{a}}}\,$ (defined in Example~\ref{ex:1.1} (e)) on the
  integrally closed integral domain $D^{\star_{a}} = D^{[\star]}$ (cf.
  Example~\ref{ex:1.1} ((c) and (d))).  Set
\[
\begin {array} {rl}
\mbox{Kr}(D,\star) := \{ f/g \;  \,|\, & f,g \in D[X] \setminus \{0\} \;\; \mbox{ \it and there
exists } \; h \in D[X] \setminus \{0\} \;   \\ &
  \mbox{ \it such that } \; (\boldsymbol{c}(f)\boldsymbol{c}(h))^\star
  \subseteq (\boldsymbol{c}(g)\boldsymbol{c}(h))^\star \,\} \, \cup\,
  \{0\}\,.
\end{array}
  \]
Then we have:

\bf (a) \it $\;$ ${\textstyle\rm Kr}(D,\star)$ is a B\'ezout domain with quotient field $K(X)\,,$
called \rm the Kronecker function ring of $D$ with respect to the semistar operation $\star\,.$
\rm

\bf (b) \it $\;$ ${\textstyle\rm Na}(D,\star) \subseteq {\textstyle\rm Kr}(D,\star)\,.$

\bf (c) \it $\;$ ${\textstyle\rm Kr}(D,\star) = {\textstyle\rm Kr}(D,\star_a) = {\textstyle\rm
Kr}(D^{\star_a},\dot{\star}_{a})\,$.

\bf (d) \it $\;$ For each $F \in \boldsymbol{f}(D)$\,:
  $$ F{\textstyle\rm Kr}(D,\star) \cap K = {\textstyle\rm Kr}(D,\star_a)
\cap K = F^{\star_{a}}\,.$$

\bf (e) \it $\;$ If $\,F := (a_{0},a_{1},\ldots, a_{n}) \in
 \boldsymbol{f}(D)$\, and $\,f(X) :=a_{0}+ a_{1}X +\ldots +a_{n}X^n \in
 K[X]\,,$ then:
 $$
 F{\textstyle\rm Kr}(D,\star) = f(X){\textstyle\rm Kr}(D,\star) =
 {\boldsymbol{c}}(f){\textstyle\rm Kr}(D,\star)\,.  $$
 \vskip -0.5cm \hfill $\square$
  \end{lemma}

The notion that we recall next is essentially due to P. Jaffard \cite{Jaffard} (cf.  also
\cite{Koch:1997}, \cite{Koch:2000}, \cite{FL2}). Let $\star$ be a semistar operation on $D$ and
let $V$ be a valuation overring of $D$.  We say that $V$ is a \it $\star$-valuation overring of
$D$ \rm if, for each $F \in \boldsymbol{f}(D)\,$, $ F^\star \subseteq FV\,$ (or equivalently,
$\star_f \leq \star_{\{V\}}$, where $\star_{\{V\}}$ is the semistar operation of finite type on
$D$ defined by:
\[ E^{\star_{\{V\}}} := EV = \cup \{FV \;| \; \; F \subseteq E,\; F \in \boldsymbol{f}(D)\}\,, \]
for each $E \in \overline{\boldsymbol{F}}(D)$; cf. Example~\ref{ex:1.1} (g) and \cite[Example 2.5
(1) and Example 3.6]{FL1}).

 \noindent  Note that a valuation overring $V$ of $D$ is a
$\star$-valuation overring of $D$ if and only if $V^{\star_{f}} =V$. (The ``only if''
part is obvious; for the ``if'' part recall that, for each $F \in \boldsymbol{f}(D)$,
there exists a nonzero element $x \in K$ such that $FV = xV$, thus $ F^\star \subseteq
(FV)^{\star_{f}} = xV^{\star_{f}}= xV = FV$.)

 We collect in the following lemma the main properties of the $\star$-valuation
overrings.

%LEMMA 1.8

\begin{lemma} \label{lm:V} \rm (\cite[Proposition 3.3, Proposition
3.4, Theorem 3.5]{FL2}) \it Let $\star$ be a semistar operation of an integral domain $D$ with
quotient field $K$ and let $V$ be a valuation overring of $D$.  Then:

\bf (a) \it $\,$ $V$ is a $\star$-valuation overring of $D$ if and only if $V$ is a
$\star_a$-valuation overring of $D$.

\bf (b) \it $\,$   $V$ is a $\star$-valuation overring of $D$ if and only if there exists a
valuation overring $W$ of ${\textstyle\rm Kr}(D,\star)$ such that $W \cap K = V\,;$ moreover, in
this case, $ W = V(X)\,.$

\bf (c) \it $\,$ $\Kr(D, \star) = \cap \{ V(X) \; | \;\; V \mbox{ is a $\ \star$--valuation
overring of } D \}\,.$

\bf (d) \it $\,$ Assume that $\, \star =\star_{a}\, $ and that $\calV$ is the set of all the
$\star$--valuation overrings of $D\,.$ For each $F \in \boldsymbol{f}(D)$,
\[ F^\star = F^{\star_{\cal V}} :=
\cap \{FV \; | \; \,V \mbox{ is a $\star$--valuation overring of $D$} \}\,,
\]
thus an e.a.b. semistar operation on $\,D\,$ is always equivalent to an a.b. semistar operation
on $\,D\,.$ \vskip -0.4cm \hfill $\square$
\end{lemma}

%LEMMA 1.9

\begin{lemma} \label{lm:M*} \rm (\cite[Theorem 4.3]{FL3}) \it Let
$\star$ be any semistar operation defined on an integral domain $D$ and let $\star_{a}$ be the
e.a.b. semistar operation of finite type associated to $\star\,.$ Assume that $\star =
\star_{f}\,.$ Then:

\bf (a) \it $\;$ $\calM(\star_a) \subseteq \{N \cap D \; | \;\; N \in \Max(\Kr(D, \star))\,\}.$

\bf (b) \it $\;$ For each $Q \in \calM(\star_a)$ there exists a $\star$--valuation overring
$\,(V, M)\,$ of $\,D\,$ such that $\,M \cap D = Q\,$ (i.e., $\,V$ dominates $D_{Q}$\,). \hfill
$\square$
\end{lemma}

Although the essential results of the theory developed in the present paper concern
finite type semistar operations, we will consider general semistar operations not only
in order to establish the results in a more general and natural setting, but also
because one the most important example of semistar operation, the (semi)star operation
$v$, is not, in general, of finite type.  The alternative use of the (semi)star
operations $v$ and $t$\,  ---\ in our case of $\star$ and $\star_f$\
---\, helps for a better understanding of the motivations and the applications of
the theory presented in this paper.

\section{Characterization of P$\star$MDs}\label{sec:2}
In this Section we prove several characterizations for an integral domain to be a
P$\star$MD, when $\,\star\,$ is a semistar operation.

We start with a first theorem in which some of the statements gene\-ra\-li\-ze some of the
classical characterizations of the P$v$MDs (cf.  M. Griffin \cite[Theorem 5]{Griffin:1967}, R.
Gilmer \cite[Theorem 2.5]{Gilmer:1970}, J. Arnold and J. Brewer \cite[Theorem
3]{Arnold/Brewer:1971}, J. Querr\'e \cite [Th\'eor\`{e}me 3, page 279]{Querre} and B.G. Kang
\cite[Theorem 3.5, Theorem 3.7]{Kang:1989}).

%THEOREM 2.1
\begin{theorem}\label{pr:3} \it Let $D$ be an integral domain and
$\star$ a semistar operation on $D$.  The following are equivalent:
\brom \rm \item \it $D$ is a P$\star$MD;

\rm \item \it $D_Q$ is a valuation domain, for each $Q\in\calM(\star_f)$;

\rm \item \it $\Na(D,\star)$ is a Pr\"{u}fer domain;

\rm \item \it $\Na(D,\star)=\Kr(D,\tilde{\star})$;

\rm \item \it $\tilde{\star}$ is an e.a.b. semistar operation;

\rm \item \it $\star_f$ is stable and e.a.b. \erom

\noindent In particular $D$ is a P$\star$MD if and only if it is a P$\tilde{\star}$MD

\end{theorem}

\begin{proof}
\bf (i) $\Rightarrow$ (ii)\rm. Let $Q\in\calM(\star_f)$ and let $J$ be a finitely
generated ideal of $D_Q$, then $J=ID_Q$ for some $I\in\boldsymbol{f}(D)$, \cite[Theorem
4.4]{GILMER:1972}.  Since $I$ is $\star_f$--invertible, then $J=ID_Q$ is invertible,
and hence principal, in the local domain $D_Q$ (Lemma~\ref{le:1} (i) $\Rightarrow$ (ii)
and \cite[Corollary 7.5]{GILMER:1972}).  As a consequence $D_Q$ is a local B\'ezout
domain, i.e., $D_Q$ is a valuation domain.

\vskip.2cm

{\bf (ii) $\Rightarrow$ (i)} is a consequence of Lemma~\ref{le:1} ((ii) $\Rightarrow$
(i)), since we are assuming that $D_Q$ is a valuation domain, for each
$Q\in\calM(\star_f)$. A direct proof is the following.  Let $ I \subseteq D $ be a
finitely generated ideal. For each $Q\in\calM(\star_f)$, we have: $$(II^{-1})D_{Q}
=(ID_{Q})(I^{-1}D_{Q}) = (ID_{Q})(ID_{Q})^{-1} =D_{Q}\, ,$$ \noindent hence $\,
II^{-1}\not\subseteq Q\,$, thus $(II^{-1})^\star= D^\star$ (Lemma~\ref{lm:1.2} (a)).

\vskip.2cm

\bf (ii) $\Rightarrow$ (iii)\rm.   The maximal ideals of the Nagata ring $\Na(D,
\star)$ are of the form $\overline{Q} := Q[X]_{N_{D}(\star)}$, for each
$Q\in\calM(\star_f)$, and we have $\Na(D,\star)_{\overline{Q}}= D_{Q}(X)$
(Lemma~\ref{lm:1.4} (f)).  If $D_{Q}$ is a valuation domain, then $D_{Q}(X)$ is also a
valuation domain and hence $\Na(D, \star)$ is a Pr\"{u}fer domain.

\vskip.2cm

\bf (iii) $\Rightarrow$ (iv)\rm.  By assumption and Lemma~\ref{lm:1.4} (i) we have that
$\Na(D, \star) = \Na(D, \tilde{\star})$ is a Pr\"{u}fer domain.  Moreover, from the
definition of $\,\tilde{\star}\,$ and from Lemma~\ref{lm:1.4} (f), (g) and (h), we
deduce that $D_{Q}$ is a $\tilde{\star}$--valuation overring of $D$, for each
$Q\in\calM(\star_f)$.  Since $\Kr(D, \tilde{\star})=\cap \{V(X) \; | \; \; V $ is a
$\tilde{\star}$--valuation overring of $D \} \,$ (Lemma \ref{lm:V} (c))\,, we obtain
that $\Kr(D, \tilde{\star}) \subseteq \cap\{D_Q(X)\; | \; \;Q\in\calM(\star_f)\}=
\Na(D, \star)\,$ (Lemma~\ref{lm:1.4} (c))\,, and thus $\Kr(D, \tilde{\star}) = \Na(D,
\star)$.

\vskip.2cm

\bf (iv) $\Rightarrow$ (v)\rm.  From the equality $\Kr(D, \tilde{\star}) = \Na(D,
\star)$ and from Lemma \ref{lm:1.4} (e) and Lemma~\ref{le:K} (d) we deduce that
$\tilde{\star} = (\tilde{\star})_{a}$.

\vskip.2cm

\bf (v) $\Rightarrow$ (ii).\rm  We recall that the following
  statements are equivalent:

{\bf (1)} $D_Q$ is a valuation domain, for each $Q\in\calM(\star_f)$;

{\bf (2)} $FD_Q$ is an invertible ideal of $D_Q$, for each $F\in\boldsymbol{f}(D)$ and for each
$Q\in\calM(\star_f)\,$;

{\bf (3)} $FD_Q$ is a \it quasi--cancellation ideal \rm of $D_Q$ (i.e., $FGD_Q\subseteq{F}HD_Q$
implies $GD_Q\subseteq{H}D_Q$ when $G$, $H\in\boldsymbol{f}(D)$), for each
$F\in\boldsymbol{f}(D)$ and for each $Q\in\calM(\star_f)$.

Note that \, (1) $\Leftrightarrow$ (2) \, since, in a local domain, for a finitely generated
ideal, invertible is equivalent to principal \cite[Corollary 7.5]{GILMER:1972}.
(2)$\Leftrightarrow$(3): this is a consequence of a                                                 %n unpublished
result by Kaplansky (cf. \cite[Exercise 7, p.  67]{GILMER:1972}, \cite[Theorem
1]{Anderson/Anderson:1984} and \cite[Theorem 13.8]{HK}).

Therefore, in order to prove that $D_Q$ is a valuation domain, for each $Q\in\calM(\star_f)$, we
show that:
\[
FGD_Q\subseteq{F}HD_Q\;\Rightarrow\;GD_Q\subseteq{H}D_Q\,,
\]
for all $F$, $G, H \in\boldsymbol{f}(D)\,.$ Now, from the assumption and from Lemma~\ref{lm:1.4}
(b), for all $E \in \boldsymbol{\overline{F}}(D)$ we have: $$ E^{\tilde{\star}} = \cap \{ ED_{Q}\;
| \;\; Q \in \calM(\star_{f}) \}\,, \;\,\mbox{and} \;\, E^{\tilde{\star}}D_{Q} = ED_{Q}\,, \;
\mbox{for each} \; Q \in \calM(\star_{f})\,.$$ \noindent  Hence, if $FGD_Q\subseteq{F}HD_Q$, then
$FG\subseteq{F}HD_Q$ and so there exists $t \in D\setminus Q$ such that $tFG \subseteq FH$. In
particular, $(FtG)^{\tilde{\star}} \subseteq (FH)^{\tilde{\star}}$, hence by assumption
$(tG)^{\tilde{\star}} \subseteq H^{\tilde{\star}}$.  From the previous remark we deduce that
$tGD_{Q}\subseteq HD_{Q}$, for each $Q \in \calM(\star_{f}) $, that is $GD_{Q}\subseteq HD_{Q}$,
because $tD_{Q }= D_{Q}$.

\vskip.2cm

\bf (vi) $\Rightarrow$ (v)\rm. Note that $\star_f$ is always quasi--spectral
(Remark~\ref{rk:1.3}) and that a semistar operation is spectral if and only if is
quasi--spectral and stable \cite[Theorem 4.12 (3)]{FH}. Therefore
\[
\star_f\mbox{ is stable }\Leftrightarrow\;\star_f\mbox{ is spectral.}
\]
Since $\,\tilde{\star} = (\star_f)_{sp}\leq\star_f\,,$ then:
\[
\star_f \mbox{ is stable }\Leftrightarrow\;\star_f=(\star_f)_{sp}=\tilde{\star}\,.
\]

\vskip.2cm

\bf (v) $\Rightarrow$ (vi)\rm.  Assume that $\,\widetilde{\star}\,$ is an e.a.b.
semistar operation (of finite type) on $\,D\,$, \, hence $\,\widetilde{\star}
=\star_{\Delta}\,$, \, where $\,\Delta:= \mathcal{M}(\star_{f}) =
\mathcal{M}(\widetilde{\star})\, $ and (by (v) $\Rightarrow$ (ii)) $D_{Q}\,$ is a
valuation domain, for each $\,Q \in \Delta\,.$
%Henceforth, $\, F^{\widetilde{\star}} = \cap \{ FD_{Q} \; | \;
%\; Q \in \Delta \}\,$, \, for each $\, F \in \boldsymbol{f}(D)\,.$

\bf Claim. \rm \it Let $\,\star\,$ be a semistar operation on an integral domain $D$.
If $\,\tilde{\star}\,$ is e.a.b., then $\,\star\,$ and $\,\star_{f}\,$ are a.b. (hence,
in particular, e.a.b.).\rm

 Assume that $\,\tilde{\star}\,$ is
e.a.b. (note that for a semistar operation of
  finite type, like $\tilde{\star}$, the notions of e.a.b. and a.b.
  coincide).  Henceforth (by (v) $\Rightarrow$ (i)) $D$ is a P$\star$MD,
  and thus each nonzero finitely generated ideal in $D$ is
  $\star_{f}$--invertible.  Let $E \in \boldsymbol{f}(D)$, and suppose
  that $(EF)^{{\star}}\subseteq (EG)^{{\star}}$, for all $F,G \in
  \boldsymbol{F}(D)$ [respectively, $F,G \in \boldsymbol{f}(D),$ for the
  e.a.b. case].  Since $E \in \boldsymbol{f}(D)$, then there exists a
  nonzero $d\in D$ such that $I:= dE$ is a nonzero finitely generated
  ideal in $D$.  Let $J\in \boldsymbol{f}(D)$ be such that $(IJ)^{\star} =
  D^\star$.
Then: $$
\begin{array}{ll}
(EF)^{{\star}}\subseteq (EG)^{{\star}} &\Rightarrow \; d(EF)^{{\star}}\subseteq d(EG)^{{\star}}\;
\Rightarrow \;
(IF)^{{\star}}\subseteq (IG)^{{\star}} \; \Rightarrow \;\\
 &\Rightarrow \; J(IF)^{{\star}}\subseteq J(IG)^{{\star}}
 \;\Rightarrow \; (J(IF)^{{\star}})^\star\subseteq (J(IG)^{{\star}})^\star
 \; \Rightarrow\\
 &\Rightarrow \; (JIF)^\star\subseteq (JIG)^\star \; \Rightarrow \;
 ((JI)^{\star}F)^\star \subseteq ((JI)^{\star}G)^\star \; \Rightarrow \\
  &\Rightarrow \; F^\star \subseteq G^\star\,.
  \end{array}
  $$

  Therefore $\,\star\,$ is  a.b..  Since $\,\widetilde{(\star_{f})}=
  \tilde{\star}\,$, from the above argument we deduce also that
  $\,\star_{f}\,$ is  a.b..

Under the present assumption, by the Claim and by \cite[Proposition 4.5 (5)]{FL1} we have that $\,
\star_{a} = (\star_{f})_{a}=\star_{f}\, $ is an a.b. semistar operation of finite type on $\,D\,$.
\, Therefore $\,\star_{f}= \star_{\mathcal{W}}\,,$ \, for some set $\,\mathcal{W}\,$ of valuation
overrings of $\,D\,$ \cite[Proposition 3.4]{FL2} (i.e. $\,F^{\star} = \cap \{ FW \; | \; \; W \in
\mathcal{W} \}\,$, \, for each $\, F \in \boldsymbol{f}(D)\,$).

Furthermore, note that, in the present situation, (by \cite[Corollary 3.8]{FL1} and (v)
$\Rightarrow$ (iv)) we have:
$$
\begin{array}{ll}
    \cap \{W(X) \; | \; \; W \in \mathcal{W} \} &= \Kr(D, \star_{a}) =\Kr(D,
    \star) \supseteq \\ &\supseteq  \Kr(D, \widetilde{\star}) =
 \Na(D, \star) = \cap \{D_{Q}(X) \;
| \; \; Q \in \Delta \}.
\end{array}
$$
\noindent Since $\, \Na(D, \star)\,$ is a Pr\"ufer domain and, by \cite[Theorem 3.9]{FL3},
$\,\mbox{\rm Max}(\Na(D, \star)$) = $ \{QD_{Q}(X)\cap \Na(D, \star) \;| \;\; Q \in \Delta\}\,$\,
then, for each $\,W \in \mathcal{W}\,$,\, there exists a prime ideal $\,Q \in \Delta\,$ and a
prime ideal $\, H\,$ in $\, D_{Q}(X)\,$, such that $\,W(X) = (D_{Q}(X))_{H}\,$.  \, Therefore, we
have that $\, W = W(X) \cap K = (D_{Q}(X))_{H} \cap K \supseteq D_{Q}\,$.  Since, for each $\,Q
\in \Delta\,,$\, $\,D_{Q}\,$ is a valuation domain, then there exists a prime ideal $\,Q'
\subseteq Q\,$ of $\, D\,$ such that $\, W = D_{Q'}\,.$ \, Set $\, \Delta' := \{Q'\; | \;\; D_{Q'}
= W \,,\; \mbox{for some} \; W \in \mathcal{W}\}\,.$\, Therefore, we have $\, \star_{\Delta}=
\widetilde{\star} \leq \star_{f}= \star_{\Delta'}\,$ \, (note that, by construction of $\,
\Delta'\,,$ \, $\, \Delta' \subseteq \Delta^{\downarrow}\,$).\, On the other hand, $\, \Delta =
\mathcal{M}(\star_{f})= \mathcal{M}(\star_{\Delta'}) \subseteq \Delta'\,$ and so
$\,{\Delta}^{\downarrow} \subseteq {\Delta'}^{\downarrow}\,.  $ \, From the previous remarks, we
deduce that $\,{\Delta}^{\downarrow} = {\Delta'}^{\downarrow}\,$ and so we conclude that
$\,\widetilde{\star} = \star_{\Delta} = \star_{\Delta'}= \star_{f} = (\star_{f})_{a}\,$.

The last statement of the theorem follows easily from the equivalence (i) $\Leftrightarrow$ (iv)
and from Lemma~\ref{lm:1.4} (i).
 \end{proof}

%REMARK
\begin{remark}\label{re:11}
As a consequence of the proof of the previous theorem, we have that:
$$
\mbox{$D$ is P$\star$MD }\;\;\; \Leftrightarrow \;\;\; \mbox{$\Na(D,
\star) = \Kr(D,\star)$}\,.
$$
As a matter of fact, when $D$ is P$\star$MD, then $\tilde{\star}=\star_{f }= (\star_{f
})_{a}= \star_{a}$ and so $\Na(D, \star) = \Kr(D,\star)$ (and conversely).
\end{remark}

Recently, W. Fanggui and R.L. McCasland \cite[Section 2]{Fanggui/McCasland:1999} have introduced,
studied and characterized the integral domains that are P$w$MD, where $\,w\,$ is the (semi)star
operation considered in Remark~\ref{rk:1.3} that, in our notation, coincides with$\,\tilde{v} \,
(= t_{sp})$.  They observed that, for a given integral domain $D$,
 $$ D \mbox{ is a P$w$MD }\;\Rightarrow \; D \mbox{ is a P$v$MD }\,.$$
The following corollary to Theorem~\ref{pr:3} shows, among other properties, that this
implication is in fact an equivalence, reobtaining a result proved by D.D. Anderson and S.J. Cook
\cite[Theorem 2.18]{Anderson/Cook:2000} that a nonzero fractional ideal is $t$--invertible if and
only if is $w$--invertible. This property was generalized in \cite[Proposition 4.25]{FH}.

 % COROLLARY TO THEOREM 2.1
\begin{corollary}\label{cor:2.2} \it Let $D$ be an integral domain. The
following are equivalent:
\brom \rm \item \it $D$ is a P$v$MD;
\rm \item \it $\Na(D, v)=\Kr(D, t_{sp})$; \rm \item \it $t_{sp}$ is an e.a.b. semistar operation.
\erom

\noindent In particular $D$ is a P$v$MD if and only if it is a P$t_{sp}$MD.

\end{corollary}
\begin{proof}  It is a straightforward consequence of the previous
theorem, after observing that $\, \tilde{v}\, = (v_{f})_{sp} = t_{sp}\,.$
\end{proof}

%REMARK 2.3
\begin{remark}\label{re:2.3}
\bf (1) \rm  Note that, if $\tilde{v} = (v_{f})_{sp} = t_{sp} =\tilde{t}\,$ is an e.a.b.
(semi)star operation on a domain $D$, then the $v$--operation is also e.a.b. operation on $D$,
but the converse is not necessarily true \cite[page 418, Theorem 34.11 and Exercise 5 page
429]{GILMER:1972}.
%%%More generally, for each semistar operation, the following holds:
%%%
%%%\it Let $\,\star\,$ be a semistar operation on an integral domain $D$.  If $\,\tilde{\star}\,$ is
%%%e.a.b. , then $\,\star\,$ is a.b. (hence, in particular, e.a.b.).\rm
%%%
%%% Assume that $\,\tilde{\star}\,$ is
%%%e.a.b., hence $D$ is a P$\star$MD (Theorem~\ref{pr:3} ((v) $\Rightarrow$ (i))),
%%% and thus
%%%each nonzero finitely generated ideal in $D$ is $\star_{f}$--invertible.  Let $E \in
%%%\boldsymbol{f}(D)$, and suppose that $(EF)^{{\star}}\subseteq (EG)^{{\star}}$, for all $F,G \in
%%%\boldsymbol{F}(D)$ [respectively, $F,G \in \boldsymbol{f}(D),$ for the e.a.b. case].  Since $E \in
%%%\boldsymbol{f}(D)$, then there exists a nonzero $d\in D$ such that $I:= dE$ is a nonzero finitely
%%%generated ideal in $D$.  Let $J\in \boldsymbol{f}(D)$ be such that $(IJ)^{\star} = D^\star$.
%%%Then: $$
%%%\begin{array}{ll}
%%%(EF)^{{\star}}\subseteq (EG)^{{\star}} &\Rightarrow \; d(EF)^{{\star}}\subseteq d(EG)^{{\star}}\;
%%%\Rightarrow \;
%%%(IF)^{{\star}}\subseteq (IG)^{{\star}} \; \Rightarrow \;\\
%%% &\Rightarrow \; J(IF)^{{\star}}\subseteq J(IG)^{{\star}}
%%% \;\Rightarrow \; (J(IF)^{{\star}})^\star\subseteq (J(IG)^{{\star}})^\star
%%% \; \Rightarrow\\
%%% &\Rightarrow \; (JIF)^\star\subseteq (JIG)^\star \; \Rightarrow \;
%%% ((JI)^{\star}F)^\star \subseteq ((JI)^{\star}G)^\star \; \Rightarrow \\
%%%  &\Rightarrow \; F^\star \subseteq G^\star\,.
%%%  \end{array}
%%%  $$
%%%We recall that for a semistar operation of finite type, like $\tilde{\star}$, the notions of
%%%e.a.b. and a.b. coincide.

\bf (2) \rm Recall that if $D$ is an integrally closed
   integral domain and if $\,D = \cap_{\alpha} V_{{\alpha }}\,$ can be
   represented as the intersection of a family of essential valuation
   overrings (e.g. if $D$ is a P$v$MD) then the a.b. (semi)star operation
   $\,\star_{\calW}\,,$ where $\calW := \{V_{\alpha}\}$ (Example~\ref{ex:1.1} (g.3)), is
   equivalent to the $\,v\,$ (semi)star operation \cite[Proposition
   44.13]{GILMER:1972}.  In particular, in a P$v$MD, $\,t_{sp}=\tilde{t}\,$
   is equivalent to $\,v\,$, i.e. $\,\tilde{t} = t\,$, since
   $\,\tilde{t}\,$ is a (semi)star operation of finite type
   (Remark~\ref{rk:1.3}).
\newline
Note that, in this context, Zafrullah \cite[Theorem
   5]{Zafrullah:1985} has proved the following general result: \it Let $D$ be an
   integral domain and $\Delta$ a set of prime ideals of $D$ such that $ D
   = \cap \{ D_{P} \, |\; P\in \Delta \}$.  Then the (semi)star operation
   $\,\star_{\Delta}\,$ is equivalent to the $\,v\,$ (semi)star operation on $D$
   if and only if, for each $F \in \boldsymbol{f}(D)$ and for each $P\in
   \Delta$, $\,FD_{P} = F^{v}D_{P}\,$.\rm \, (It is obvious that when $D_{P}$
   is a valuation domain, then $\,FD_{P} = F^{v}D_{P} = (FD_{P})^{v}\,$,
   because $\,FD_{P}\,$ is a principal ideal in  $D_{P}$.)

\bf (3) \rm For Pr\"{u}fer domains, J. Arnold \cite[Theorem 4]{Arnold:1969} has proved that,
   if $D$ is an integral domain, then:
 $$
 \begin{array}{ll}
     D \mbox{ is a Pr\"{u}fer domain } &\Leftrightarrow\; \Na(D, d) =
 D(X) \mbox{ is a Pr\"{u}fer domain }\; \Leftrightarrow\\
 &\Leftrightarrow\; \Na(D, d)
 =\Kr(D, b)\,.
 \end {array}$$
 Note that the previous equivalence follows from Theorem~\ref{pr:3} ((i)
 $\Leftrightarrow$ (iii) $\Leftrightarrow$ (iv)), since if $D$ is
 Pr\"{u}fer then $\,d = \tilde{d} = b\,$ and if $\Na(D, d) =\Kr(D, b)$ then
$\,d = \tilde{d} = b_{a} = b\,$ is an e.a.b. (semi)star operation.

\end{remark}

Next result gives a positive answer to the problem of the
 ``ascent'' of the P$\star$MD property.
%PROPOSITION 2.4

\begin{proposition}\label{pr:5}
Let $\star$ be a semistar operation defined on an integral domain $D$ and let $T$ be an overring
of $D$.  Denote simply by $\dot{\star}$ the semistar operation ${\dot{\star}}^{\mbox{{\tiny
$T$}}}$ on $T\,$ (Example~\ref{ex:1.1} (e)).\, Assume that $D$ is a P$\star$MD, then $T$ is a
P$\dot{\star}$MD.
 \end{proposition}
\begin{proof}
To avoid the trivial case, we can assume that $\,T\,$ is different from the quotient field of
$D\,.$  Let $H$ be a prime ideal of $T$ which is a maximal element in the set of nonzero ideals
of $T$ with the property that $H^{\dot{\star}_f}\cap{T}=H$, i.e. $H$ is a
quasi--${\dot{\star}_f}$--maximal of $T$.  We want to show that $T_H$ is a valuation domain
(Theorem~\ref{pr:3} ((ii)$\Rightarrow$(i))).  If we consider the prime ideal $\,Q:=H\cap{D}\,$ of
$\, D\,$, then $\, Q \,$ is nonzero, since $\,{D}_Q\subseteq{T}_H\,,$ and moreover:
\[
\begin{array}{ll}
Q^{\star_f}\cap{D} =(H\cap{D})^{\star_f}\cap{D}
&\subseteq{H}^{\star_f}\cap{D} =\\
  & = H^{\dot{\star}_f}\cap{T}\cap{D}
= H\cap{D} =Q \subseteq{Q}^{\star_f}\cap{D}\,,
\end{array}
\]
and thus $Q$ is a prime quasi--$\star_f$--ideal of $D$. If $Q$ is not a quasi--$\star_f$--maximal,
then there exists a prime ideal $P$ such that $Q\subseteq{P}$ and $P = P^{\star_f}\cap{D}$
(Lemma~\ref{lm:1.2} (a)).  Now we have:
\[
D_P\subseteq{D}_Q\subseteq{T}_H
\]
with $D_P$ valuation domain, because $D$ is a P$\star$MD (Theorem~\ref{pr:3}
((i)$\Rightarrow$(ii))).  We conclude immediately that $T_H$ is a valuation domain.
\end{proof}

%COROLLARY 2.5
\begin{corollary}\label{co:7}
Let $\star$ be a semistar operation defined on an integral domain $D$. Assume that $D$ is a
P$\star$MD and denote simply by  ${\dot{\star}}$ the (semi)star operation
$\dot{\star}^{\mbox{{\tiny D$^{\star}$}}}$ on $D^\star$ (Example~\ref{ex:1.1} (e)).  Then
$D^\star$ is a P$\dot{\star}$MD.
\end{corollary}
\begin{proof}
The statement is a straightforward consequence of Proposition \ref{pr:5} (taking $T$ =
$D^{\star}$).
\end{proof}

Next goal is to study the ``descent'' of the P$\star$MD property. The following lemma is required
in the proof of next proposition.

%LEMMA 2.6
\begin{lemma}\label{lm:2.4} Let $T$ be an overring of an integral domain
$D$ and let $\star$ be a semistar operation on $T$.  The semistar operations of finite type \rm \
(\d{$\star$})$_{f}$ \ \it and \ \rm \d{($\star_{f}$)} \ \it (both defined on $D$) coincide. (For
the sake of simplicity, we will simply denote  by \ \d{$\star$}$_{f}$ \ this semistar operation.)
\end{lemma}
\begin{proof} Let $E \in \boldsymbol{\overline{F}}(D)$, then
 $$
 \begin{array}{ll}
 E^{\mbox{\tiny (\d{$\star$})$_{f}$}}
 &=\cup\{F^{\mbox{\tiny \d{$\star$}}} \;| \; \, F\subseteq E, \; F \in
    \boldsymbol{f}(D) \} = \cup \{(FT)^{\star} \;|\; \, F\subseteq E, \; F
    \in \boldsymbol{f}(D) \} \subseteq \\
 &\subseteq \cup \{H^{\star} \;|\; \, H\subseteq ET, \; H \in \boldsymbol{f}(T) \} \ (=(ET)^{\star_{f}} = E^{ \mbox{
 \tiny \d{($\star_{f}$)} }}) = \\
 &= \cup \{(FT)^{\star} \;|\;  \, F\subseteq ET, \; F \in
 \boldsymbol{f}(D) \} \subseteq \\
&\subseteq \cup \{(FT)^{\star} \;|\; \, F\subseteq E, \; F \in \boldsymbol{f}(D) \} =
E^{\mbox{\tiny(\d{$\star$})$_{f}$}}\,,
\end{array}
$$
\noindent since, if $\, F \subseteq ET\, $ with $F \in \boldsymbol{f}(D) $,  it is possible to
find $\,E_{0}\subseteq E \, $ with $E_{0} \in \boldsymbol{f}(D)$ and $\, F \subseteq E_{0}T \,$,
\, therefore $\, (FT)^{\star} \subseteq (E_{0}T)^{\star}\, $.
\end{proof}

%PROPOSITION 2.7
\begin{proposition}\label{pr:8}
Let $T$ be a flat overring of an integral domain $D$. Let $\star$ be a semistar operation on
$T$.  Assume that $T$ is a P$\star$MD. Denote simply by \d{$\star$} the semistar operation
\d{$\star$}$_{\mbox{{\tiny D}}}$ on $D$ (Example~\ref{ex:1.1} (e)).  Then $D$ is a P\
\hskip-0.09cm\d{$\star$}MD\,.
\end{proposition}
\begin{proof}
Let $Q\in\calM({\mbox{\d{$\star$}}}_f)$, then by Lemma~\ref{lm:2.4} we have
$\,Q^{{\mbox{\d{$\star$}}}_f}\cap{D} $ = $(QT)^{\star_f}\cap{D}=Q \,$.  In particular
$QT\neq{T}$, hence there exists $H\in\calM(\star_f)$ such that $H\supseteq{QT}$ and so
$H\cap{D}\supseteq{Q}$.  Note that $\, (H\cap{D})^{{\mbox{\d{$\star$}}}_f} =
((H\cap{D})T)^{\star_f}\,,$ and since $H\in\calM(\star_f)\,,$ then:
\[
H\cap{D} \subseteq ((H\cap{D})T)^{\star_f}\cap{D} \subseteq H^{\star_f}\cap{D}
=H^{\star_f}\cap{T}\cap{D} =H\cap{D}\,.
\]
Henceforth, $H\cap{D}$ is a quasi--${\mbox{\d{$\star$}}}_f$--prime of $D$ and so $H\cap D = Q\,.$
Therefore, we conclude that $\calM({\mbox{\d{$\star$}}}_f)$ coincides with the contraction to $D$
of the set $\calM(\star_f)$.  Since $T_H$ is a valuation domain, for each $H\in\calM(\star_f)\,,$
and $T$ is $D$--flat then, by \cite[Theorem 2]{R}, we conclude that $D_{H\cap{D}}=T_H$ is also a
valuation domain, and so $D$ is a P${\mbox{\d{$\star$}}}$MD (Theorem~\ref{pr:3} ((ii)
$\Rightarrow$ (i))).
\end{proof}

%REMARK 2.8
\begin{remark}\label{re:9}
 \textbf{(1)}
 Note that, in Proposition~\ref{pr:8}, the hypothesis that $T$ is
$D$--flat is essential (cf.  also \cite[Theorem 27.2]{HK}). For example, let $(T,M)$ be
a discrete 1-dimensional valuation domain with residue field $k$.  Let $k_0$ be a
proper subfield of $k$ and assume that $k$ is a finite field extension of $k_0\,.$ Set
 \[
 \begin{xy}
\xymatrix{
 D:=\varphi^{-1}(k_0)\ar@{->>}[rr]\ar@{^(->}[dd]&&k_0\ar@{^(->}[dd]\\
 \\
  T\ar@{->>}[rr]^\varphi&&T/M=k
}\end{xy}
\]
\noindent Then $D$ and $T$ are local with the same maximal ideal $M\,,$ which is a finitely
generated ideal both in $D$ and in $T\,,$ \cite[Theorem 2.3]{F}.  Let $\star := b \ (= d)$ be the
identical (semi)star operation on the valuation domain $T$.  Then \d{$\star$}$_{\mbox{{\tiny D}}}
=\star_{\{T\}}$, i.e. $E^{\mbox{\d{$\star$}}_{D}} =ET\,,$ for each
$E\in\overline{\boldsymbol{F}}(D)$.  Obviously $T$ is a (local) Pr\"{u}fer domain, but $D$ is not a
P$\star_{\{T\}}$MD, since $M\in \calM(({\mbox{\d{$\star$}}}_{\mbox{{\tiny D}}})_f) =
\calM((\star_{\{T\}})_f)$ but $D_M=D$ is not a valuation domain.

\bf{(2)} \rm Note that, from Proposition~\ref{pr:5} and
 Example~\ref{ex:1.1} (e.6), if $\star$ is a semistar operation on the overring $T$ of $D\,,$
 if $\mbox{\d{$\star$}} = {\mbox{\d{$\star$}}}_{\mbox{{\tiny D}}}$ and if
 $D$ is a P\d{$\star$}MD, then
 $T$ is a P$\star$MD.
\end{remark}

%EXAMPLE 2.9
\begin{example} \label{ex:10} \it When $\,\star\,$ is a semistar operation,
  a P$\star$MD, is not ne\-ces\-sa\-ri\-ly integrally closed.  \rm (Note
  that if $\,\star\,$ is a semistar operation on an integral domain
  $\,D\,$ and $\,D\,$ is a P$\star$MD, then $\,D^{\star}\,$ must be
  integrally closed by Corollary 2.3 and \cite[Theorem 34.6, Proposition
  34.7 and Theorem 34.11]{GILMER:1972}; in particular, if $\,\star\,$ is a
  star operation on $\,D\,,$ then $\,D\,$ is integrally closed.)
\newline
Let $\,D\,$ be a non integrally closed integral domain and let
 $\,\Delta\,$ be a nonempty finite set of nonzero prime ideals of $\,D\,$
 with the following properties:

\bf (a) \rm $\,D_P\,$ is a valuation domain, for each $P \in \Delta\,$;

\bf (b) \rm $\,D_{P'}\,$ and $\,D_{P''}\,$ are incomparable, if $P' \neq P''\; \mbox{ and }\;
P',  P''\in \Delta\,.$
\newline
Let $ \star := \star_{{\Delta}}$ be the spectral semistar operation on $\,D\,$ associated to
${\Delta}$ (Example~\ref{ex:1.1} (f)).  Since $\, D^{\star}= \cap\{D_{P}\; | \;\; P \in
\Delta\}\, ,\,$ Max$(D) = \{PD_{P}\cap D^{\star}\; | \;\; P \in \Delta\}\,$ and $\, D^{\star}\,$
is a semilocal B\'ezout domain \cite[Theorem 107]{Kaplansky:1970}, then clearly $\, D \subsetneq
D^{\star}\,$ and $\,D^{\star}\,$ is flat over $\,D\,$ \cite[Theorem 2]{R}.

Let $\boldsymbol{\ast} := \dot{\star} = \dot{\star}^{{\mbox{\tiny \it D}}^{\star}}$ denote the
(semi)star operation defined on $\, D^{\star}\,$ induced by $\, \star\,$ (Example~\ref{ex:1.1}
(e)), then $\, D^{\star}\,$ is trivially a P$\boldsymbol{\ast}$MD, since $\, D^{\star}\,$ is a
B\'ezout domain.  Denote simply by \d{$\boldsymbol{\ast}$} the semistar operation
\d{${\boldsymbol{\ast}}$}$_{\mbox{\tiny \it \tiny D}}\,$ on $\, D\,$ induced by $\,
\boldsymbol{\ast}\,$ (Example~\ref{ex:1.1} (e)) then, by Proposition~\ref{pr:5}, $\,D\,$ is a
P\d{$\boldsymbol{\ast}$}MD, but by assumption is not integrally closed.  \ Note that it is easy
to verify that, in the present situation, $\, \star = \mbox{\d{$\boldsymbol{\ast}$}}$
  since, for each $E \in \boldsymbol{\overline{F}}(D)$, we have: $$
  E^{\mbox{\scriptsize \d{$\boldsymbol{\ast}$}}} =
  (ED^{\star})^{\mbox{\scriptsize $\boldsymbol{\ast}$}}
= (ED^{\star})^{\star} = (ED)^{\star} = E^{\star}\,.$$ Therefore $\,D\,$ is a P$\star$MD but, by
assumption, it is not integrally closed. In particular $D$ is not a P$v$MD.
\end{example}

The following explicit construction produces an example similar to the situation described in
previous Remark~\ref{re:9} (1).

 %EXAMPLE 2.10
\begin{example}\label{ex:2.7} Let $K$ be a field and $X$, $Y$ indeterminates over
$K$.  Set $F:=K(X)$ and $D:=K+YF[[Y]]$.  It is well known that $D$ in an integrally closed
1-dimensional non-valuation local domain with maximal ideal $M:=YF[[Y]]$ and that
$V:=K[X]_{(X)}+YF[[Y]]$  is a 2-dimensional valuation overring of $D$ with maximal ideal $N:=
XK[X]_{(X)}+YF[[Y]]$, \cite[Section 17, Exercises 11, 12, 13, 14 and page 231]{GILMER:1972}.
Note that $M$ is also an ideal inside $V$, and precisely $M$ is the height 1 prime ideal of $V$.
\newline
Consider the semistar operation $\star:=\star_{\{V\}}$ on $D$ (cf. Example~\ref{ex:1.1} (g)).  It
is of finite type and induces over $\,V=D^\star\,$ the identity (semi)star operation $\,d_{V}\,$
on $\,V$, i.e. $ \dot{\star}\, ( = {\dot{\star}}^{{\mbox{\tiny \it V}}}) = d_{V}\,.$ Henceforth
$D^\star$ is a P${\dot{\star}}$MD, in fact it is a valuation domain.
\newline
Note that $D$ is not P$\star$MD, because the only maximal (quasi)$\star$--ideal is $M$, since
$M^\star = MV = M$, and because $D= D_{M}$ is not a valuation domain (Theorem~\ref{pr:3}, (i)
$\Leftrightarrow$ (ii)).
\newline
Keeping in mind Proposition~\ref{pr:8}, note also that $V$ is not $D$--flat by \cite[Theorem
2]{R}, because it is easy to see that $V_{M} = F[[Y]] \supsetneq D_{M} =D$. Moreover, if
$\,\delta:= d_{V}\,$ is the identical (semi)star operation on $\,V\,$,  then the semistar
o\-pe\-ra\-tion $\,$ \d{$\delta$} :=\d{$\delta$}$_{_{\mbox{\tiny \it \tiny D}}}$ $\,$ on $\,D\,$
induced by $\,\delta\,,$ defined in Example~\ref{ex:1.1} (e), coincides with $\,\star\,.$
\newline
Note, also, that in the present situation $\, \tilde{\star} = \star_ {sp}\,,$ since $\, \star =
\star_f\,;\,$ moreover $\, \star_ {sp} = d_{D}\, $ the identical (semi)star operation on $D$,
since $\,{ \cal{M}}(\star_f) = \{M\}\, $ and $D_M = D$.  Furthermore, $\, \tilde{\star} =d_{D}
\,$ is not an e.a.b. (semi)star operation on $D$ (cf.  also Theorem~\ref{pr:3} ((i)
$\Leftrightarrow$ (v)), because of the equivalence $(1) \, \Leftrightarrow \, (3)$ in the proof
(v) $\Rightarrow$ (ii) of Theorem~\ref{pr:3} and because $D = D_M$ is not a valuation domain.
\end{example}

The previous example shows that if $D^\star$ is a P$\dot{\star}$MD then $D$ is not
ne\-ces\-sa\-rily a P$\star$MD. This fact induces to strengthen the condition ``$D^\star$ is
P$\dot{\star}$MD'' for characterizing $D$ as a P$\star$MD and it suggests (in the finite type
case) the use of the semistar operation $\star_{sp}$ (or, equivalently, $\tilde{\star}$) instead
of $\star$.

%PROPOSITION 2.11
\begin{proposition}\label{co:11} Let $\star$ be a semistar operation defined
    on an integral domain $D$.
With the notation of Lemma~\ref{lm:1.4}, we have: $$ D \, \mbox{ \it is a P$\star$MD} \; \;
\Leftrightarrow \; \; D \, \mbox{ \it is a P$\tilde{\star}$MD} \; \; \Leftrightarrow \; \;
D^{\tilde{\star}}\, \mbox{ \it is a P$\dot{\tilde{\star}}$MD}\,. $$
\end{proposition}
\begin{proof}
>From Theorem~\ref{pr:3} and Corollary~\ref{co:7}  we deduce immediately that:
 $$ D\, \mbox{ is a P$\star$MD} \; \Leftrightarrow \; D\, \mbox{ is a P$\tilde{\star}$MD} \; \Rightarrow \;
 D^{\tilde{\star}}\, \mbox{ is a P$\dot{\tilde{\star}}$MD}\,. $$
\noindent Set $\,\widetilde{D} := D^{\tilde{\star}}\,$. By Lemma~\ref{lm:1.4} (h) we know that
$$\calM(\dot{\tilde{\star}}) = \{\widetilde{Q} := QD_{Q}\cap \widetilde{D} \; | \;\; Q \in
\calM(\star_{f})\}\; \, \mbox{and} \; \, \widetilde{D}_{\widetilde{Q}} = D_{Q}\,,$$ for each $\ Q
\in \calM(\star_{f})\,.$ \, Assume that $\widetilde{D}\, \mbox{ is a P$\dot{\tilde{\star}}$MD}\,,$
then ${\widetilde{D}}_{\widetilde{Q}} = D_{Q}\,$ is a valuation domain, for each $Q \in
\calM(\star_{f})\,$, by Theorem~\ref{pr:3} ((i) $\Rightarrow$ (ii))) applied to
$\,\widetilde{D}\,.$ We conclude that $D$ is a P$\tilde{\star}$MD from Theorem~\ref{pr:3} ((ii)
$\Rightarrow$ (i)) and from Lemma~\ref{lm:1.4} (g).
\end{proof}

Next example shows that the flatness hypothesis in Proposition~\ref{pr:8} is essential also
outside of a pullback setting (cf. for instance Remark~\ref{re:9} (1) and Example~\ref{ex:2.7}).

%EXAMPLE 2.12
\begin{example}\label{ex:2.9} Let $\,T\,$ be an overring of an integral
domain $\,D\,$ and let $\, \star := \star_{\{T\}}\,$ be the semistar operation of finite type on
$\,D\,,$ defined in Example~\ref{ex:1.1} (g) with $\, \mathcal{T} := \{ T \}\,.$ Assume that
$\,T\,$ is integral over $\,D\,$ and that $\,D \neq T\,,$ then $\,D\,$ is not a P$\star$MD even
if $\,T\,$ is a Pr\"{u}fer domain.

 Note that, as in Example~\ref{ex:2.7}, if $\,\delta:= d_{T}\,$ is the
 identical (semi)star o\-pe\-ra\-tion on $\,T\,$,   then the semistar operation
 $\,$ \d{$\delta$} :=\d{$\delta$}$_{_{\mbox{\tiny \it \tiny D}}}$ $\,$ on
 $\,D\,$ induced by $\,\delta\,,$ defined in Example~\ref{ex:1.1} (e), coincides
 with $\,\star\,.  \,$ Moreover, since $\,T\,$ is integral over $\,D\,$
 then, by the lying-over theorem, we have: $$\mbox{Max}\{P \in
 \mbox{Spec}(D) \; | \;\; 0 \neq P \, \mbox{ and }
 \, PT \cap D \neq D \} = \mbox{Max}(D)\,.$$
\noindent
 Therefore, by \cite[Chapitre II, \S 3, N.  3, Corollaire 4]{Bourbaki},
 Lemma 1.2 (c) and Remark 1.3, we have:
 $$ \tilde{\star} = \star_{sp} = d\,,$$
 \noindent where $\,d := d_{D}\,$ is the identical (semi)star operation on $\,D\,,$ and so $\,
D^{\tilde{\star}}= D\,.  \,$ By using Proposition 2.6, we have: $$ D \, \mbox{ is a P$\star$MD
}\;\;\Leftrightarrow \;\; D \, \mbox{ is a P$d$MD $\,$ (i.e.  $\,D\,$ is a Pr\"{u}fer domain) }\,,$$
and this is excluded if $\, D \neq T\,.$

More generally, the previous argument shows that:

\it Let $\,T\,$ be a proper integral overring of an integral domain $\,D\,.\, $ Assume that there
exists a semistar operation $\,\ast\,$ on $\,T\,$ such that $\,T\,$ is a P$\ast$MD. Then $\,D\,$
is not a P${\star}$MD, for any semistar operation $\,\star\,$ on $\,D\,$ such that $\, \star \leq
\star_{\{T\}}\, ( \, \leq \mbox{\d{$\ast$}}\ )\,.$ \rm

In fact, recall that if $\star_{1}$ and $\star_{2}$ are two semistar operations on an integral
domain $D$, if $\, \star_{1} \leq \star_{2}\,$ and if $D$ is a P$\star_{1}$MD, then $D$ is also a
P$\star_{2}$MD.
 Therefore, it is
sufficient to show that $\,D\,$ is not a P${\star_{\{T\}}}$MD and this fact follows from the
equivalence proved above, since $\,D\,$ is not a Pr\"{u}fer domain because, by assumption, $\,T \neq
D\,$ is integral over $\,D\,.$
\end{example}

In case of star operations, next goal is to characterize P$\star$MDs in terms of P$v$MDs.  We
start with few general remarks concerning the ``star setting''.

%REMARK 2.13
\begin{remark}\label{qu:14} \rm Let $\star,\star_{1}\,$ and $\,
\star_{2}$ be \sl star operations \rm on an integral domain $D\,.$ We denote by Spec$_{\star}(D)$
the set of all prime ideals $P$ of $D$, such that $P^\star = P$, then obviously: $$ \star_{1}
\leq \star_{2} \;\; \Rightarrow \;\; \mbox{Spec}_{\star_{2}}(D) \subseteq
\mbox{Spec}_{\star_{1}}(D)\,.$$ A prime ideal $P$ of an integral domain $D$ is called \it a
valued prime
\rm if $D_{P}$ is a valuation domain. \\
Let $\star$ be a star operation on $D$ and assume that $D$ is a P$\star$MD (hence, in particular,
a P$v$MD).  Then, by \cite[Proposition 4.1]{Mott/Zafrullah:1981}, a prime ideal of $D$ is valued
if and only if it is $t$--ideal.  As a consequence, under the present assumptions, the valued
prime ideals of $D$ are inside $\Spec_{t}(D)$.  Moreover, since $\star_{f} \leq t\,$, for each
star operation $\, \star\, $ on $\, D\,$ \cite[Theorem 34.1 (4)]{GILMER:1972}, then
$\Spec_{t}(D)\subseteq \Spec_{\star_f}(D)$.  On the other hand, as $D$ is P$\star$MD, then each
maximal $\star_f$--ideal is valued, hence $\calM(\star_f) \subseteq \Spec_{t}(D)$, but this means
$\Spec_{\star_f}(D)\subseteq \Spec_{t}(D)$, which implies that $\Spec_{\star_f}(D) =
\Spec_{t}(D)$.

\end{remark}

In the following proposition we prove that the implication (i) $\Rightarrow$ (ii), due
to Kang \cite[Theorem 3.5]{Kang:1989}, can be inverted, obtaining a new
characterization of a P$\star$MD which is related to \cite[Proposition
21]{Garcia/Jara/Santos:1999} (cf. also \cite[Theorem 17.1 ii)]{HK}):

%PROPOSITION 2.14
 \begin{proposition} \label{pr:3.4} Let $\star$ be a star
operation on an integral domain $D\,.$ The following statements are equivalent: \brom \rm \item
\it $\, D\, $ is a P$\star$MD. \rm \item \it $\, D\, $ is a P$v$MD and $\,\tilde{\star} = t\,.$
\rm \item \it $\, D\, $ is a P$v$MD and $\,{\star}_{f} = t\,.$ \erom
    \end{proposition}
    \begin{proof} (i) $\Rightarrow$ (ii).  Since $\,\Spec_{\star_f}(D) =
    \Spec_{t}(D)\,$ (Remark~\ref{qu:14}), then $\, \tilde{\star} =
    \tilde{t}\,.$ Moreover a P$\star$MD is a P$v$MD and, in a P$v$MD, $\,
    \tilde{t} =t\,$ \cite[Theorem 3.5]{Kang:1989}.

    (ii) $\Leftrightarrow$ (iii).  It is a consequence of Remark~\ref{re:11}.

    (iii) $\Rightarrow$ (i) Since $\, v_{f} = t = \star_{f}\,$, then the
    conclusion follows immediately from the fact that
    the notions of P$\star$MD and P$\star_{f}$MD
    coincide, for each semistar operation $\star$.
    \end{proof}

>From the previous result it is possible to find star operations $\,\star\,$ on an integral domain
$D$ such that $D$ is a P$v$MD, but $D$ is not a P$\star$MD. For instance, if $D$ is a Krull non
Dedekind domain, then obviously $D$ is a P$v$MD but not a P$\star$MD, if $\,\star\,$ coincides
with $\,d\,$ the identical star operation on $D$, since $d = d_{sp}$ and, in a Krull domain $D$,
$t = d$ if and only if $D$ is a Dedekind domain, \cite[Theorem 34.12 and Theorem
43.16]{GILMER:1972}.

Next example describes a more general situation.

\begin{example}
Let $K$ be a field and $X$ and $Y$ be indeterminates over $K$. Let us consider two distinct
maximal ideals $M_1$ and $M_2$ of $K[X,Y]$. Let $S:=K[X,Y]\setminus(M_1\cup{M_2})$ be a
multiplicative closed subset of $K[X,Y]$ and let $D:=S^{-1}K[X,Y]$. Thus $D$ is a Noetherian
Krull domain, hence $D$ is a P$v$MD. Moreover $D$ is semilocal with maximal ideals
$N_1=S^{-1}M_1$ and $N_2=S^{-1}M_2$ (note that $D_{N_1}$ and $D_{N_2}$ are not valuation domains).
\newline
Let us consider the spectral star operation $\star$ on $D$ defined by the subset
$\Delta:=\Spec(D)\setminus\{N_2\}$, i.e. $\star=\star_\Delta$ as in Example~\ref{ex:1.1} (f). It
is not difficult to show that $\star\neq{d}$ (in fact $(N_2)^d=N_2\neq{D}=N_2^\star$) and $D$ is
not P$\star$MD (as $N_1$ is a maximal $\star$--ideal and $D_{N_1}$ is not a valuation domain).
\end{example}

Next result makes more precise the statement of Proposition~\ref{pr:3.4}
    in case $\star$ coincides with the identical star operation (cf.
    also \cite[Theorem 17.3]{HK}).

    %PROPOSITION 2.15
    \begin{proposition} \label{cor:3.5} Let $D$ be an integral domain, then the
    following are equivalent:
     \brom \rm
    \item \it $\, D\, $ is a Pr\"{u}fer domain.  \rm
    \item \it $\, D\, $ is integrally closed and $\,d = t\,.$ \rm
    \item \it $\, D\, $ is integrally  closed and has a unique star operation
    of finite type. \erom
    \end{proposition}
    \begin{proof} It is obvious that (i) $\Rightarrow$ (ii) $\Leftrightarrow$
    (iii), since it is well known that for each star operation of finite
    type $\,\star\,$ of an integral domain, $\, d \leq \star \leq t\,$,
    \cite[Theorem 34.1 (4)]{GILMER:1972}.  Finally (ii) $\Rightarrow$ (i)
    because, under the present assumptions, for each nonzero ideal $I$ of $D$, we
    have:
    $$ I_{t} = I = \cap \{ID_{M}\; | \;\, M \in \mbox{Max}(D)\}$$
    \noindent where, obviously, Max$(D) ={\calM}(d)={\calM}(t)\,,$ and thus
    the conclusion follows from \cite[Theorem 3.5]{Kang:1989}.
    \end{proof}

    %REMARK 2.16
      \begin{remark}\label{rk:2.16}
      \bf (1) \rm From the previous result we deduce that, \it in a Pr\"{u}fer
      domain, any two star operations are equivalent (in fact, both are
      equivalent to $\,v\,$ \rm (Proposition~\ref{cor:3.5} ((i) $\Rightarrow$
      (ii)))  \it  and each star operation $\,\star\,$ is a.b. (in fact,
      $\,\tilde{\star}\,$ is a.b. \rm  (cf.  Remark~\ref{re:2.3})),
      \cite[Proposition 32.18]{GILMER:1972}.
\newline
      The last part of the statement follows from the fact that each
      localization of $D$ is a valuation domain, thus the star operation
      $\,\tilde{\star}\,$ is necessarily a $\,\star_{\calW}$--operation, for
      some family $\, \calW\,$ of valuation overrings of $D$ (Example~\ref{ex:1.1} (g.3)).
\newline
    In relation with the first part of the statement note that,
    for each star ope\-ra\-tion $\,\star\,$ on a Pr\"{u}fer domain, we have
    $\,\star_{f} = t = d = b\,$ and thus $\,
    \tilde{\star}=\widetilde{(\star_{f})} =$ $= \tilde{d} = d = b = \star_{f} = t
    \,$.

    \bf (2) \rm Note that the statement in (1) is not a characterization of
    Pr\"{u}fer domains, since there exists an integrally closed non--Pr\"{u}fer
    integral domain such that any two star operations are equivalent and
    each star operation is a.b. \cite[Section 32, Exercise 12]{GILMER:1972}
    and \cite[Proposition 24]{Okabe/Matsuda:1997}.
\newline
    On the other hand, for an integral domain $D$, we have:

\vskip.2cm
 \centerline{\it $D$ is Pr\"{u}fer if and only if each semistar operation on $D$ is a.b.\rm}

\vskip.2cm \noindent     By an argument as in (1), we have that if $D$ is Pr\"{u}fer
    then each semistar ope\-ra\-tion on $D$ is a.b. Conversely, for each
    prime ideal $P$ of $D$, if $\,\star_{\{D_{P}\}}\,$ is an a.b. operation
    then, by the equivalence (1) $\Leftrightarrow$ (3) in the proof of
    Theorem~\ref{pr:3} ((v) $\Rightarrow$ (ii)), we deduce that $D_{P}$ is a
    valuation domain.
    \end{remark}

The following remark provides a ``quantitative
    information'' about the size of the set of all the semistar operations
    $\star$ on a given integral domain $D$ for which $D$ is a P$\star$MD.

    %REMARK 2.17
    \begin{remark}\label{qu:15}
Let ${\boldsymbol{\cal P}}(D)$ be the set of all semistar operations of finite type on $D$ such
that $D$ is a P$\star$MD and let ${\boldsymbol{\cal B}}(\Spec(D))$ be the set of all the subsets
of $\Spec(D)\,.$ Then, the map:
 $$ \mu:{\boldsymbol{\cal P}}(D) \rightarrow{\boldsymbol{\calB}}(\Spec(D))\,,
 \; \; \star \mapsto \calM(\star_{f})\,,$$
defines a surjection onto the set ${\boldsymbol{\cal M}}(D) \,(\subseteq {\boldsymbol{\cal
B}}(\Spec(D)))\,$ of all the subsets of $\Spec(D)$ that are quasi--compact and that are formed by
valued incomparable prime ideals of $D$ \cite[Corollary 4.6]{FH}.  Obviously $\mu(\star_{1}) =
\mu(\star_{2})$ if and only if ${\tilde{\star}}_{1} = {\tilde{\star}}_{2}\,$
(Remark~\ref{rk:1.3}).  Note that the map:
 $$\mu': {\boldsymbol{\cal M}}(D) \rightarrow
 {\boldsymbol{\cal P}}(D)\,, \;\; \calM \mapsto \star_{\calM}\,,$$
is such that $\mu \circ \mu'$ is the identity.
\end{remark}

Next goal is to give a characterization of a P$\star$MD,
 when $\star$ is a semistar operation, in terms of polynomials, by
 generalizing the classical cha\-racte\-ri\-za\-tion of Pr\"{u}fer domain in
 terms of polynomials given by R. Gilmer and J. Hoffman \cite[Theorem
 2]{Gilmer/Hoffmann:1974}.  Note that similar properties, in the ``star
 setting'', were already considered by J. Mott and M. Zafrullah
 \cite[Theorem 3.4]{Mott/Zafrullah:1981} and by E. Houston, S.J. Malik
 and J. Mott \cite[Theorem 1.1]{Houston/Malik/Mott:1984}.
\newline
 Let $D$ be an integral domain with quotient field $K\,,$
 recall that an ideal $I$ of a polynomial ring $D[X]$ is called \it an
 upper to $0$ in $D[X]$ \rm if there exists a nontrivial ideal $J$ in
 $K[X]$ such that $J \cap D[X] = I$.  Note that a nontrivial primary
 ideal $H$ of $D[X]$ is an upper to $0$ if and only if $H\cap D = 0$.

%THEOREM 2.18

\begin{theorem}\label{ht:11} Let $\star$ be a semistar operation
defined on an integral domain $D$ with quotient field $K$.  The following statements are
equivalent:
\brom
 \rm \item \it $D$ is a P$\star$MD;
 \rm \item \it $D^{\tilde{\star}}$ is integrally closed
 (i.e. $D^{\tilde{\star}} = D^{[\tilde{\star}]}$)
  and, for each $\,I\,$ upper to
 $\,0\,$ in $\,D[X]\,,$ we have $\,I\Na(D,\star)=\Na(D,\star)\,$ (or,
 equivalently, there exists $\,f\in{I}\,$ such that
 $\,\boldsymbol{c}(f)^\star = D^\star$);
 \rm \item \it $D^{\tilde{\star}}$ is integrally closed
 (i.e. $D^{\tilde{\star}} = D^{[\tilde{\star}]}$) and, for each nonzero prime ideal $\,H\,$ of
 $\,D[X]\,$ such that $\,H\cap{D}=0\,$, we have $\,H\Na(D,\star) =
 \Na(D,\star)\,$ (or, equivalently, there exists $\,f\in{H}\,$ such that
 $\,\boldsymbol{c}(f)^\star= D^\star$);
 \rm \item \it $D^{\tilde{\star}}$ is integrally closed
 (i.e. $D^{\tilde{\star}} = D^{[\tilde{\star}]}$) and, for all nonzero elements $\,a$, $b\in{D}\,,$
 the prime ideal $\,H:=(aX+b)K[X]\cap{D}[X]\,$ of $\,D[X]\,$ is such that
 $\,H\Na(D,\star)=\Na(D,\star)\,$ (or, equivalently, there exists
 $\,f\in{H}\,$ such that $\,\boldsymbol{c}(f)^\star=D^\star$).  \erom
\end{theorem}

\begin{proof}
(i) $\Rightarrow$ (ii).  We know, from Corollary~\ref{co:7}, that $\,D^\star\,$ is a
P$\dot{\star}$MD, where $\,\dot{\star} = \dot{\star}^{{\mbox{\tiny D}}^\star}\,$ defines a star
operation on $\,D^\star\,$ (when restricted to $\boldsymbol{F}(D^\star))$, and hence
$\,D^\star\,$ is integrally closed \cite[Corollary 32.8 and Theorem 34.11]{GILMER:1972}.  The
same argument can be applied to $\,\tilde{\star}\,$ and $\,D^{\tilde{\star}}\,$.  Moreover, since
$\,\tilde{\star}\,$ is stable by Example~\ref{ex:1.1} (f.2), then  we deduce that
$D^{\tilde{\star}} = D^{[\tilde{\star}]}$ or, equivalently, that $D^{\tilde{\star}}$ is
integrally closed (Example~\ref{ex:1.1} (c.2)).

\noindent
 From Theorem~\ref{pr:3} ((i) $\Rightarrow$ (iv)) we know that, in the
 present si\-tua\-tion $\,\Na(D,\star)$ $=\Kr(D,\tilde{\star})\,$.  Therefore, if $\,I
 := hK[X] \cap D[X]\,,$ with $\,h\,$ a non constant polynomial of
 $\,K[X]\,,$ by Lemma~\ref{le:K} (e), we have:
 $$
  \begin{array}{ll}
 I\Na(D, \star) &=\; I\Kr(D, \tilde{\star}) \supseteq
 \{\boldsymbol{c}(f)\Kr(D, \tilde{\star})
 \; | \;\; f \in I \} =\\
 & =\; \{\boldsymbol{c}(hg)\Kr(D, \tilde{\star}) \; | \;\; g \in K[X]\,,\;
 hg \in D[X] \} = \\
 &=\; \{\boldsymbol{c}(hg)\Na(D, \star) \; | \;\; g \in K[X]\,, \; hg \in D[X]
 \}\,.
 \end{array}
 $$
 Since $\,D\,$ is a P$\star$MD, then there exists a finitely generated
 (fractional) ideal $\,L\,$ of $\,D\,$ such that
 $\,(\boldsymbol{c}(h)L)^\star=D^\star\,$ and $\,L \subseteq (D:_{\mbox{\tiny \it K}}\boldsymbol{c}(h))\,.$

\noindent Let $\,\ell \in{K}[X]\,$ be such that $\,\boldsymbol{c}(\ell)=L\,$. Then, by the content
formula \cite[Theorem 28.1]{GILMER:1972}, for some $m\geq0\,,$ we have
\[
 \boldsymbol{c}(h)\boldsymbol{c}(\ell)\boldsymbol{c}(h)^m=\boldsymbol{c}(h\ell)\boldsymbol{c}(h)^m
\]
and so
\[
 (\boldsymbol{c}(h)\boldsymbol{c}(\ell)\boldsymbol{c}(h)^{m}L^{m})^{\star} =
  (\boldsymbol{c}(h\ell)\boldsymbol{c}(h)^mL^m)^{\star} \,.
\]
Therefore:
\[
 D^\star=(\boldsymbol{c}(h)\boldsymbol{c}(\ell))^\star=\boldsymbol{c}(h\ell)^\star\,.
\]
Set $\,f:=h\ell\,,$ since  $\,L \subseteq
 (D:_{\mbox{\tiny \it K}} \boldsymbol{c}(h))\,$ then $\, f \in
 I \subseteq I\Na(D, \star)\,.$ By the fact that
 $\, \boldsymbol{c}(f)^\star = D^\star\,$, we deduce
 $\,f\,\Na(D^\star,\star)=\Na(D^\star,\star)\,,$ and thus
 $\,I\Na(D^\star,\star)=\Na(D^\star,\star)$.

(ii) $\Rightarrow$ (iii) $\Rightarrow$ (iv) are trivial.

(iv) $\Rightarrow$ (i).  Set $\widetilde{D} := D^{\tilde{\star}}\,$ and, for each $\, Q \in
\calM(\star_{f})\,$, $\, \widetilde{Q} := QD_{Q}\cap \widetilde{D}\,.$ Note that, by
Lemma~\ref{lm:1.4} (h), $\, {{\widetilde{D}}_{\widetilde{Q}}} = D_{Q}\,.$ \, By assumption,
$\,{\widetilde{D}}\,$ (and so ${\,\widetilde{D}}_{\widetilde{Q}}\,$) is integrally closed, for
each $\, Q \in \calM(\star_{f})\,$.  In order to conclude we want to show that
$\,{{\widetilde{D}}_{\widetilde{Q}}}\,$ is a valuation domain, for each $\, Q \in
\calM(\star_{f})\,$ (Proposition~\ref{co:11} and Theorem~\ref{pr:3} ((ii) $\Rightarrow$ (i))).
Let $t:=a/b\in{K}$ with $a$, $b\in{D}$, $b\neq 0$, and let $H:=(bX-a)K[X]\cap{D}[X]\,.$ By
assumption, there exists a polynomial $\, f \in H \subseteq D[X]\,$ such that $\,
{\boldsymbol{c}}(f)^\star = D^\star\,.$ In particular we have that $\,{\boldsymbol{c}}(f) \in D
\setminus Q\,,$ since $\,Q \in \calM(\star_{f})$.  Henceforth $ f\, \in D[X] \setminus QD[X]
\subseteq \widetilde{D}[X] \setminus \widetilde{Q}\widetilde{D}[X]\,.$ Since $\,f \in H\,$ then
$f(t) = 0\,,$ this implies that $\,t\,$ or $\,t^{-1}$ is in ${{\widetilde{D}}_{\widetilde{Q}}}\,$
\cite[Lemma 19.14]{GILMER:1972}.
\end{proof}
\section{Passing through field extensions}\label{sec:3}

 In this section we deal with the preservation of the
 P$\star$MD property by ``ascent'' and ``descent'', in case of field
 extensions.  Our purpose is to generalize to the P$\star$MD case the
 following classical results concerning Pr\"{u}fer domains (cf.
 \cite[Theorem 22.4 and Theorem 22.3]{GILMER:1972}):

 \bf (1) \rm Let $D$ be an integrally closed domain with
 quotient field $K$ which is a subring of an integral domain $T$.
 Assume that $T$ is integral over $D$ and that $T$ is a Pr\"{u}fer domain,
 then $D$ is also a Pr\"{u}fer domain.

 \bf (2) \rm Let $D$ be a Pr\"{u}fer domain with quotient
 field $K$ and let $L$ be an algebraic field extension of $K$.  Then the
 integral closure $T$ of $D$ in $L$ is a Pr\"{u}fer domain.

 When we study the ``descent'' of the P$\star$MD property,
 we have to consider also a ``natural restriction'' of the semistar
 operation $\,\star\,$.  Recall that, in 1936 W. Krull \cite[Satz
 9]{Krull:1936} proved that if $D$ in an integrally closed integral
 domain with quotient field $K$, if $L$ is an algebraic field extension
 of $K$ and if $T$ is the integral closure of $D$ in $L$, then, for each
 nonzero fractional ideal $E$ of $D$,
$$ (ET)^v \cap K = E^v\,,$$ (cf.  also \cite[Lemma 3.7]{Koch:2000}). The same formula holds, when
$X$ is indeterminate over $K$ and $T:= D[X]$ (cf. \cite[Section 34, Exercise 16]{GILMER:1972}).
\newline
The following result shows that, when we assume for the ``natural restriction'' that a property
of the previous type holds, then we have a ``descent'' theorem for P$\star$MDs:

%PROPOSITION 3.1
\begin{proposition}\label{pr:12}
Let $K\subseteq{L}$ be any field extension and let $T$ be an integral domain with quotient field
$L$. Assume that $D:=T\cap{K}\neq{K}$, that $T$ is integral over $D$ and that $\star$ is semistar
operation on $T$ such that $T$ is a P$\star$MD. Define ${\underline{\star}}_{D} {\colon}
\boldsymbol{\overline{F}}(D) \to \boldsymbol{\overline{F}}(D)$ in the following way:
\[
 E^{{\underline{\star}}_{D}}:=(ET)^\star\cap{K}\,.
\]
Then:

\bf {(1)} \it the operation ${\underline{\star}}_{D}$ is a semistar operation on $D$;

\bf {(2)} \it $D$ is a P${\underline{\star}}_{D}$MD.
\end{proposition}

\begin{proof}   (1) It is obvious that, if $E,
F \in \boldsymbol{\overline{F}}(D)\,,$ then $E\subseteq{F}$ implies
$E^{{\underline{\star}}_{D}}\subseteq{F}^{{\underline{\star}}_{D}}$.  Moreover, if $E\in
\boldsymbol{\overline{F}}(D)$ and $x\in{K}$, $x\neq 0$, then:
\[
\begin{array}{ll}
 (E^{{\underline{\star}}_{D}})^{{\underline{\star}}_{D}}
 &= \, (((ET)^\star\cap{K})T)^\star\cap{K}
 \subseteq(((ET)^\star)T)^\star\cap{K} = \\
 &= \, (ET)^{\star\star}\cap{K}
 =(ET)^\star\cap{K}
 =E^{{\underline{\star}}_{D}}\,;\\
 (xE)^{{\underline{\star}}_{D}}
 &= \, (xET)^\star\cap{K}
 =x(ET)^\star\cap{K}
 =x((ET)^\star\cap{K}) =\\
 &= xE^{{\underline{\star}}_{D}}\,.
\end{array}
\]

 (2) Since $T$ is a P$\star$MD, then $T_H$ is a valuation domain, for each
 $H\in\calM(\star_f)\,$ (Theorem~\ref{pr:3} ((i) $ \Rightarrow $ (ii))). By the
 assumption that $D\subseteq{T}$ is an integral extension, we know that,
 if we denote by $P$ the prime ideal $H\cap{D}$, then $D_P=T_H\cap{K}$
 and so $D_P$ is a valuation domain \cite[Theorem 22.4, Proposition
 12.7]{GILMER:1972}.  To conclude we need to show that, for all
 $Q\in\calM((\underline{\star}_D)_f)\,,$ there exists
 $H\in\calM(\star_f)$ such that $H\cap{D}=Q$.

 \bf{Claim 1}.  \it For each $H\in\calM(\star_f)\,,$ the prime ideal $P:=H\cap{D}$ has the
following property:
\[
 P^{({\underline{\star}}_D)_f}\cap{D}=P\,.
\]
\rm As a matter of fact, since $D=T\cap{K}\,,$ then $D^{{\underline{\star}}_D}=T^\star\cap{D}$
and so:
\[
 P^{({{\underline{\star}}_D})_f} =(PT)^{\star_f}\cap{K}
 =(((H\cap{D})T)^{\star_f})\cap{K} \subseteq{H}^{\star_f}\cap{K}\,.
\]
Therefore:
\[
 P^{({{\underline{\star}}_D})_f}\cap{D} \subseteq{H}^{\star_f}\cap{D}
 =H^{\star_f}\cap{T}\cap{D} =H\cap{D} =P\,.
\]

\bf{Claim 2}.  \it If $Q\in\calM(({\underline{\star}}_D)_f)\,,$ then there exists
$H\in\calM(\star_f)$ such that $Q\subseteq{H}\cap{D}$.

\rm Since $(QT)^{\star_f}\cap{K}=Q$, then also $(QT)^{\star_f}\cap{D}=Q$. Take
$L:=(QT)^{\star_f}\cap{T}\,,$ then it is easy to see that $L^{\star_f}\cap{T}=L$ and
$L\cap{D}=Q$.  Therefore, by Lemma~\ref{lm:1.2} (a), $L$ is contained in some $H\in\calM(\star_f)$
with $H\cap{D}\supseteq{L}\cap{D}=Q$.  \newline From the previous claims, we deduce that
$\calM(({\underline{\star}}_D)_f)$ coincides with the contraction of $\calM(\star_f)$ into $D$.
This is enough to conclude.
\end{proof}

>From Krull's result concerning the $\,v\,$ (semi)star operation cited before
Proposition~\ref{pr:12}, we deduce immediately that:

%COROLLARY 3.2
\begin{corollary}\label{cor:3.2}
Let $K\subseteq{L}$ be an algebraic field extension, let $T$ be an integral domain with quotient
field $L$, set $D:=T\cap{K}$.  Assume that  $T$ is the integral closure of $D$ in $L$ and that
$T$ is a P$v$MD. Then $D$ is a P$v$MD. \hfill $\square$
\end{corollary}

In \cite[Section 11]{Prufer} H. Pr\"{u}fer showed that the integral closure of a Pr\"{u}fer domain
[respectively, a P$v$MD] in an algebraic field extension is still a Pr\"{u}fer domain [respectively,
a P$v$MD].  An explicit proof of a stronger form of this result with different techniques was
given recently by F. Lucius \cite[Theorem 4.6 and Theorem 4.4]{Lucius:1998} (cf. also
\cite[Theorem 3.6]{Koch:2000}).  A generalization to the case of P$\star$MDs, when $\star$ is a
semistar operation, is proven next.

%THEOREM 3.3
\begin{theorem}\label{pr:13}
Let $K\subseteq{L}$ be an algebraic field extension. Let $D$ be an integral domain with quotient
field $K$. Assume that $\star$ is a semistar operation on $D$ such that $D$ is a P$\star$MD. Let
$T$ be the integral closure of
 $D^\star$ into $L$, and let:
\[
\calW:=\{W\mbox{ is valuation domain of }L\; | \,\, W\cap{K}=D_Q\mbox{ , for some }
Q\in\calM(\star_f)\}.
\]
 For each
$E\in\boldsymbol{\overline{F}}(T)$, set:
\[
 E^{{\overline{\star}}^T}:=\cap\{EW\;|\,\,W\in\calW\}\,.
\]

\bf{(1)} \it The operation $\, \overline{\star}^T \, $ is an a.b. (semi)star operation on $T$.

\bf{(2)} \it $T$ is a P${\overline{\star}}^T$MD.
\end{theorem}

\begin{proof}
 First at all, note that $\calW$ is nonempty (Theorem~\ref{pr:3} ((i)
 $\Rightarrow$ (ii)) and \cite[Theorem 20.1 ]{GILMER:1972}),  $\, T=\cap\{W \; | \;\; W\in\calW\}\,$
 (\cite[Theorem 19.6 and Theorem 19.8]{GILMER:1972}) and $
 T^{{\overline{\star}}^T} = T\,.$

(1)  is a straightforward
 consequence of Example~\ref{ex:1.1} (a) and (g.3).

(2)  It is sufficient to show that $\Na(T,{\overline{\star}}^T)$ is Pr\"{u}fer domain
(Theorem~\ref{pr:3} ((iii) $\Rightarrow$ (i))).  Since $D$ is a P$\star$MD, then $\Na(D,\star)$
is a Pr\"{u}fer domain (Theorem~\ref{pr:3} ((i) $\Rightarrow$ (iii))), so it is the same for its
integral closure $\overline{\Na(D,\star)}$ in the algebraic field extension $L(X)$ of $K(X)\,,$
\cite[Theorem 22.3]{GILMER:1972}.  If we show that $\overline{\Na(D,\star)}\subseteq
\Na(T,{\overline{\star}}^T)$, then we conclude \cite[Theorem 26.1 (1)]{GILMER:1972}.  In order to
prove this fact it is enough to note that:

(a) \hspace*{10pt}$\Na(T,{\overline{\star}}^T)$ is integrally closed in $L(X)\,$;

(b) \hspace*{10pt}$\Na(D,\star)\subseteq\Na(T,{\overline{\star}}^T)\,$.

For (a), we have that
\[
 \Na(T,{\overline{\star}}^T)=\cap\{T_H(X)\; | \;\; H\in\calM(({\overline{\star}}^T)_f)\}\,, \;
 \mbox{\, (Lemma~\ref{lm:1.4} (c)),}
\]
and $T_H$ is integrally closed (and, thus, $T_{H}(X)$ is integrally closed), for each $H$, since
$T$ is integrally closed.

For (b), let $\, z=f/g\in\Na(D,\star) = \Na(D,\tilde{\star})\,,$ with $\,f, g \in{D}[X]\,$ and
$\, D^{\tilde{\star}}= {\boldsymbol{c}(g)}^{\tilde{\star}} =\cap\{ \boldsymbol{c}(g)D_Q \; | \;\;
Q\in\calM(\star_f)\}\,,$ (Lemma~\ref{lm:1.4} ((b) and (i)). Then:
\[
 T = {D}^{{\overline{\star}}^{T}} \subseteq
 {(D^{\tilde{\star}})}^{{\overline{\star}}^{T}}=
 ({{\boldsymbol{c}(g)}}^{\tilde{\star}})^{\overline{\star}^T}= \cap \{
 \boldsymbol{c}(g)W \; | \;\; W\in\calW\} = {{\boldsymbol{c}(g)}}^{{\overline{\star}}^{T}} \,,
\]
and so $1\in {\boldsymbol{c}(g)}^{{\overline{\star}}^T}$, i.e.
${\boldsymbol{c}(g)}^{{\overline{\star}}^T} = T^{{\overline{\star}}^T}=T$.  Therefore
$f/g\in\Na(T,{{\overline{\star}}^T})$.
\end{proof}

>From the previous result and Corollary~\ref{cor:3.2} we reobtain the following result (cf.  for
instance \cite[Theorem 4.6]{Lucius:1998}):

%COROLLARY 3.4
\begin{corollary}\label{cor:3.4}
Let $K\subseteq{L}$ be an algebraic field extension, let $T$ be an integral domain with quotient
field $L$, set $D:=T\cap{K}$.  Assume that $D$ is integrally closed and that $T$ is the integral
closure of $D$ in $L$. Then
 $D$ is a P$v$MD if and only if $T$ is a P$v$MD.
\end{corollary}
\begin{proof}  With the notation of
Theorem~\ref{pr:13}, it is sufficient to remark that if $v_{\mbox{\tiny $D$}}\,$ [respectively,
$\,v_{\mbox{\tiny $T$}}\,$] is the $\,v$--operation on $D$ [respectively, on $T$] and if $D$ is a
P$v_{\mbox{\tiny$D$}}$MD, then the a.b. semistar operation ${\overline{v_{\mbox{\tiny
$D$}}}}^{T}$ on the integrally closed domain $T$ is equivalent to $ v_{{\mbox{\tiny $T$}}}\,$
(Remark~\ref{re:2.3} (2)).
\end{proof}

\end{document}